\documentclass[preprint]{elsarticle}

\usepackage{amssymb,latexsym,amsmath, amsfonts}
\usepackage{amsthm,indentfirst}

\newcommand{\z}[1]{\textbf{#1}}

\newtheorem{thm}{Theorem}[section]
\newtheorem{lem}[thm]{Lemma}
\newtheorem{cor}[thm]{Corollary}
\newtheorem{rem}[thm]{Remark}

\newtheorem{dfn}[thm]{Definition}

\def\cC{{\mathfrak{C}}}
\def\N{{\mathbb{N}}}

\def\Rl{{\mathbb{R}}}
\def\Cx{{\mathbb{C}}}
\def\sH{{\mathcal{H}}}
\def\sI{{\mathcal{I}}}

\def\aM{{\mathcal{M}}}

\def\mL{{\mathcal{L}}}

\journal{Journal of Functional Analysis}

\begin{document}

\begin{frontmatter}

\title{Singular traces and perturbation formulae of higher order}

\author{Denis Potapov}
\ead{d.potapov@unsw.edu.au}
\author{Fedor Sukochev}
\ead{f.sukochev@unsw.edu.au}
\author{Alexandr Usachev\corref{cor}}
\ead{a.usachev@unsw.edu.au}
\author{Dmitriy Zanin\fnref{l1}}
\ead{d.zanin@unsw.edu.au}

\cortext[cor]{Corresponding Author}

\fntext[l1]{The research of all four authors was supported by the Australian Research Council.}

\address{School of Mathematics and Statistics, University of New South Wales, Sydney, 2052, Australia.}

\begin{abstract}
Let $H, V$ be self-adjoint operators
such that $V$ belongs to the weak trace class ideal. We prove higher order perturbation formula
$$\tau\big(f(H+V)-\sum_{j=0}^{n-1}\frac{1}{j!}\frac{d^j}{dt^j} f(H+tV)\bigg|_{t=0}\big)=\int_{\Rl} f^{(n)}(t)\,dm_n(t),$$
 where $\tau$ is a trace on the weak trace class ideal and $m_n$ is a finite measure that is not necessarily absolutely continuous. 
This result extends the first and second order perturbation formulas of Dykema and Shripka, who generalised the Krein and Koplienko trace formulas to the weak trace class ideal. 
We also establish the perturbation formulae when the perturbation $V$ belongs to the quasi-Banach ideal weak-$L_n$ for any $n \geq 1$.
\end{abstract}


\end{frontmatter}

\section{Introduction}
For a pair of self-adjoint operators $H$ and $V$ acting on a separable Hilbert space $\sH$
such that $V$ belongs to the trace class ideal $\mL_1$, Krein proved that
there is a unique function $\xi_{H,V} \in L_1 (\mathbb R)$,
called the \emph{spectral shift function,} such that the following trace formula holds:
$${\rm Tr} \left( f(H+V) - f(H)\right) = \int_{\Rl} f' (t)\, \xi (t)\, dt,$$
 for every function $f \in C^1(\Rl)$ whose derivative $f'$ admits the representation
$$f'(\lambda)=\int_{\Rl} e^{-\lambda t} d\mu(t)$$
for some finite (complex) Borel measure $\mu$ on $\Rl$
(see~\cite{Krein})  and {\rm Tr} is the standard trace on trace class operators. 
Krein's formula originated from research in physics \cite{Lifsic}. It has been applied in perturbation theory \cite{Yafaev}
and in noncommutative geometry in the study of spectral
flow \cite{ACDS, ACS}. Dykema and Skripka \cite{DS14} recently extended Krein's formula to perturbation by a weak trace class operator, where the formula now involves a trace on the weak trace class operators. We recall all traces on the weak trace class are singular, that is they vanish on trace class operators and particularly on those of finite rank. As noted in \cite{DS14} the use of traces on the weak ideal introduces new perturbation formulas
that calculate the difference up to trace class perturbation and also introduces spectral measures that are not absolutely continuous. Singular traces are important in classical and noncommutative geometry as well as in applications to physics (see, e.g., \cite{CS, C_book, LSZ} and references cited therein),
and perturbation formulas for singular traces on the weak trace class ideal open new ground for applications. In particular, 
a very recent publication \cite{PSTZ} employs singular traces of a Taylor expansion (as in Theorem~\ref{M2} below) for a concrete function $f(t)=t^p$, which are used as a technical tool for studying Fr\'echet differentiability of the $L_p$-norm of Haagerup $L_p$ spaces.

Krein's trace formula was extended to Hilbert-Schmidt operators by Koplienko~\cite{Koplienko}.
He showed that if a perturbation~$V$ is from $\mL_2$ and if~$f$ is any rational function with non-real poles which
is bounded at infinity, then the difference
$$f(H+V) - f(H) - \frac{d}{dt} \left( f (H + t V) \right)\Bigr|_{t = 0}$$
belongs to $\mL_1$ and there exists a
function~$\eta = \eta_{H, V} \in L_1({\mathbb R})$ such that
\begin{equation}
\label{Koplienko}
\mathrm{Tr} \left( f(H+V) - f(H) - \frac{d}{dt} \left( f (H + t V) \right)\Bigr|_{t = 0} \right) =
\int_{\mathbb R} f''(t) \, \eta(t)\, dt.
\end{equation}

Recently Potapov, Skripka and Sukochev~\cite{PSS} have
extended the trace formula further to the case of the $n^{\mathrm{th}}$ Schatten-von Neumann ideal $\mL_n$.
More precisely, they proved that for~$V \in \mL_n$ there exists a
function $\eta_n=\eta_{n,H,V}$ in $L_1({\mathbb R})$ such that
\begin{align}
\label{PSS}
\mathrm{Tr}\left(f(H)-\sum_{j=0}^{n-1}\frac{1}{j!}\frac{d^j}{dt^j} f(H+tV)\bigg|_{t=0}\right)=\int_{\mathbb
R} f^{(n)}(t)\eta_n(t)\,dt.
\end{align}

The main aim of the present paper is to  extend (\ref{PSS}) 
to the weak Schatten-von Neumann ideals $\mathcal{L}_{n,\infty}$, which are the $n^{\mathrm{th}}$ order convexification of the weak trace class ideal 
$\mathcal{L}_{1,\infty}$ as the Schatten-von Neumann ideals $\mathcal{L}_n$ are the  $n^{\mathrm{th}}$ order convexification of the $\mathcal{L}_1$.
The weak trace class ideal possesses a infinite lattice of traces which are well studied~\cite{CS, LSZ, SSUZ}.  
The objective to extend (\ref{PSS}) naturally therefore involves traces on the weak trace class ideal 
(which are all singular traces, that is traces vanishing on finite rank operators) instead of the classical trace.

The first step in this direction was done recently by Dykema and Skripka (see~\cite{DS14}). They proved
the analogues of Krein's and Koplienko's results for a specific class of Lorentz ideals $\aM_\psi$
and Dixmier traces, which are singular traces of a special type (see e.g.~\cite{LSZ}).
We state their results for the special case when the ideal  $\aM_\psi$ coincides with the classical Dixmier-Macaev ideal 
$\mathcal M_{1,\infty}$ of all compact operators $A \in B(\sH)$ such that 
$$
\|A\|_{\aM_{1,\infty}}:=\sup_{n \ge 1}  \frac{1}{\log(1+n)} \sum_{k=0}^{n-1}\mu(k,A) < \infty,
$$
where $\{\mu(n,A)\}_{n\ge0}$ is the sequence of singular values of a compact operator $A \in B(\sH)$.
Let $\aM^{(n)}_{1,\infty}$ denote the $n^{\mathrm{th}}$ root (or $n$-convexification) of $\aM_{1,\infty}$, in other words
$\aM^{(n)}_{1,\infty}$ consists of all operators $A$ such that $A^n \in \aM_{1,\infty}$ (see e.g.~\cite{P_env, CRSS, LSZ} for details).

\begin{thm}[{\cite[Theorems 3.3, 5.2]{DS14}}]
1. If $H$ and $V$ are self-adjoint operators and if $V$ belongs to the Lorentz ideal $\aM_{1,\infty}$,
then for every bounded trace $\tau$ on $\aM_{1,\infty}$ there exists a unique finite Borel measure~$m_1$ on $\Rl$
depending only on $H,V$ and $\tau$ such that
    \begin{equation}
      \tau \left( f(H+V) - f(H) \right) = \int_{\Rl}f' (t)\,  dm_1(t),
    \end{equation}
    for every $f \in C^3(\Rl)$.

2. If $H$ and $V$ are self-adjoint operators and if $V$ belongs to the Lorentz ideal $\aM^{(2)}_{1,\infty}$,
then for every bounded trace $\tau$ on $\aM_{1,\infty}$ there exists a unique finite Borel measure~$m_2$ on $\Rl$
depending only on $H,V$ and $\tau$ such that
    \begin{equation}
      \tau \left( f(H+V) - f(H) - \frac{d}{dt} \left( f (H + t V) \right)\Bigr|_{t = 0} \right) = \int_{\Rl}f'' (t)\,  dm_2(t),
    \end{equation}
    for every functions $f$ from
    $${\rm span} \{ \lambda \mapsto (z-\lambda)^{-k} : k\in \N, {\rm Im} (z) <0 \}.$$
\end{thm}

The result of~\cite{PSS} extends the results of Krein and Koplienko, and, in a similar fashion,
we extend here the result of Dykema and Skripka~\cite{DS14} (in the setting of self-adjoint operators) as follows:
\begin{thm}\label{M1}
Let $H, V\in B(\sH)$ be self-adjoint operators
such that~$V \in \aM^{(n)}_{1,\infty}$.
For every bounded trace $\tau$ on $\aM_{1,\infty}$ there exists a unique totally finite Radon signed measure~$m_n$ on   $\Rl$
depending only on $n,H,V$ and $\tau$ such that
\begin{align}
\tau\left(f(H+V)-\sum_{j=0}^{n-1}\frac{1}{j!}\frac{d^j}{dt^j} f(H+tV)\bigg|_{t=0}\right)=\int_{\mathbb
R} f^{(n)}(t)\,dm_n(t),
\end{align}
for every Schwartz function $f$. 
Moreover, the total variation of $m_n$ is bounded as follows: $\|m_n\| \le \|V\|^n_{\aM^{(n)}_{1,\infty}}$.
\end{thm}

To treat the case of an unbounded operator $H$ we narrow the class of functions $f$ to that considered in~\cite{Koplienko}.

\begin{thm}\label{M3}
Let $H, V$ be self-adjoint operators
such that~$V \in \aM^{(n)}_{1,\infty}$.
For every bounded trace $\tau$ on $\aM_{1,\infty}$ there exists a unique totally finite Radon signed measure~$m_n$ on $\Rl$
depending only on $n,H,V$ and $\tau$ such that
\begin{align}
\tau\left(f(H+V)-\sum_{j=0}^{n-1}\frac{1}{j!}\frac{d^j}{dt^j} f(H+tV)\bigg|_{t=0}\right)=\int_{\mathbb
R} f^{(n)}(t)\,dm_n(t),
\end{align}
for every rational function $f$ with non-real poles which
is bounded at infinity. 
Moreover, the total variation of $m_n$ is bounded as follows: $\|m_n\| \le \|V\|^n_{\aM^{(n)}_{1,\infty}}$.
\end{thm}


The following theorems provide the trace formulae for weak-$\mL_n$ ideals,
that are the proper sub-ideals in $\aM_{1,\infty}^{(n)}$ of the form
$$
\mathcal L_{n,\infty}=\mL^{(n)}_{1,\infty}:= \left\{A \in \aM_{1,\infty}^{(n)}:  \|A\|_{\mL_{n,\infty}}:=\sup_{k \ge 0}  (k+1) \mu^n(k,A) < \infty\right\}.
$$
\begin{thm}\label{M2}
Let $H, V\in B(\sH)$ be self-adjoint operators
such that~$V \in \mL_{n,\infty}$.
For every bounded trace $\tau$ on $\mL_{1,\infty}$ there exists a unique totally finite Radon signed measure~$m_n$  on $\Rl$
depending only on $n,H,V$ and $\tau$ such that
\begin{align}
\tau\left(f(H+V)-\sum_{j=0}^{n-1}\frac{1}{j!}\frac{d^j}{dt^j} f(H+tV)\bigg|_{t=0}\right)=\int_{\mathbb
R} f^{(n)}(t)\,dm_n(t),
\end{align}
for every Schwartz function $f$.
Moreover, the total variation of $m_n$ is bounded as follows: $\|m_n\| \le \|V\|^n_{\mL_{n,\infty}}$.
\end{thm}


\begin{thm}\label{M4}
Let $H, V$ be self-adjoint operators
such that~$V \in \mL_{n,\infty}$.
For every bounded trace $\tau$ on $\mL_{1,\infty}$ there exists a unique totally finite Radon signed measure~$m_n$  on $\Rl$
depending only on $n,H,V$ and $\tau$ such that
\begin{align}
\tau\left(f(H+V)-\sum_{j=0}^{n-1}\frac{1}{j!}\frac{d^j}{dt^j} f(H+tV)\bigg|_{t=0}\right)=\int_{\mathbb
R} f^{(n)}(t)\,dm_n(t),
\end{align}
for every rational function $f$ with non-real poles which
is bounded at infinity.
Moreover, the total variation of $m_n$ is bounded as follows: $\|m_n\| \le \|V\|^n_{\mL_{n,\infty}}$.
\end{thm}

\begin{rem}
Although the results of Theorems~\ref{M1} and~\ref{M3} (and also that of Theorems~\ref{M2} and~\ref{M4})
present somewhat similar formulae, the bounded operator cases turned out to be more difficult, since a broader class of functions $f$ is considered there.
In order to treat these cases we extend the theory of multiple operator integrals to the quasi-Banach setting,
using the integration techniques originated in the papers of Turpin and Waelbroeck~\cite{TW}, and Kalton~\cite{K86b} (see also~\cite{HS}).
\end{rem}

\begin{rem}\label{remark}
At a glance Theorems~\ref{M2} and~\ref{M4} are a particular case of Theorem~\ref{M1} and~\ref{M3}, since $\mL_{1,\infty}^{(n)}=\mL_{n,\infty} \subset \aM_{1,\infty}^{(n)}$.
There is a difference however. It lies in the distinction between the sets of traces on $\aM_{1,\infty}$ and $\mL_{1,\infty}$.
Indeed, if we apply Theorem~\ref{M1} to $V \in \mL_{n,\infty}$ we obtain the result of Theorem~\ref{M2}
for every bounded trace on $\mL_{1,\infty}$  that is a restriction of a bounded trace on $\aM_{1,\infty}$.
However, it is proved in~\cite[Theorem 4.7]{SSUZ}  that not every bounded trace on $\mL_{1,\infty}$ is the restriction of a bounded trace on $\aM_{1,\infty}$.

To sum up, Theorem~\ref{M2} (respectively, Theorem~\ref{M4}) provides a result which is more general  than Theorem~\ref{M1} (respectively, Theorem~\ref{M3}) applied to $\mL_{1,\infty}$.
\end{rem}
\begin{rem}
 In contrast to Theorem 1.1 from~\cite{PSS},
 we do not prove the absolute continuity of measure $m_n$ for Theorem~\ref{M2},
 since in general this is not the case (see~\cite[Proposition 4.2]{DS14}).
\end{rem}

\section{Preliminaries}
Let $C^n$ denote the space of all $n$ times
continuously differentiable complex-valued functions on $\Rl$ equipped with the usual norm and let $C_b^n$ denote the subclass of $C^n$ of bounded functions.
Also by $C(\Rl)$ we denote the space of continuous real-valued functions equipped with the standard norm.

Let $B(\sH)$ be the algebra of all bounded linear operators on a separable Hilbert space
$\sH$ equipped with the operator norm $\|\cdot\|_\infty$.

\begin{dfn}
 A trace on an ideal $\sI$ of $B(\sH)$ is a linear functional $ \tau : \sI \to \Cx$ such
that
$$\tau(AB) = \tau(BA)$$
for all $A \in \sI$ and $B\in B(\sH)$. A trace $\tau$ is called singular if it vanishes on finite rank operators.
\end{dfn}
Note, that we do not require a trace to be positive.

By ${\rm Tr}$ we denote the standard (normal) trace on $B(\sH)$.
Let $\mL_p:=\mL_p(\sH)$ denote the Schatten-von Neumann ideal, that is the set of all $A \in B(\sH)$ such that ${\rm Tr} (|A|^p) <\infty$, $0<p<\infty$.
(see, e.g., \cite{ACDS, PiXu03} and references cited therein
for basic definitions and facts).
Whereas there are no non-trivial bounded traces on $\mL_p$, $0<p<1$, 
the stock of traces on $\mL_{1,\infty}$ is plentiful.
In particular, it contains Dixmier traces~\cite{D} (see also~\cite[Section IV.2.$\beta$]{C_book}), which are bounded traces of a special form.

We now define the derivatives appearing in the perturbation formulae.

\begin{dfn}\label{G}
 Let $g: \Rl \to \Rl$. Let $\mathcal {SA}$ be the subspace of all self-adjoint operators from $B(\sH)$ and $A, B\in \mathcal {SA}$.
 Consider the function $t\mapsto g(A+tB)$ and define its derivative as follows
 $$\frac {d}{dt} \left[ g\left( A +tB \right) \right] |_{t = 0} :=\lim_{t\to0} \frac{g(A+tB)-g(A)}{t}$$
 provided the limit exists. 
 By $\frac {d^k}{dt^k} \left[ g\left( A +tB \right) \right]|_{t = 0}$, $k\ge2$ we understand the derivative
 of the function
 $$t \mapsto \frac {d^{k-1}}{dt^{k-1}} \left[ g\left(A +tB \right) \right]|_{t = 0}.$$
\end{dfn}

\begin{dfn}\label{F}
 Let $g: \mathcal {SA} \to \mathcal {SA}$ and $A\in \mathcal {SA}$.
 The function $g$ is called $k$ times Fr\'echet differentiable at $A$
 if there exist bounded linear operators
$$F_A^{(j)} : \underbrace{\mathcal {SA} \times \cdots \times \mathcal {SA}}_{j-\text{times}} \to \mathcal {SA}, \quad 1\le j \le k$$ such that
$$\lim_{B\to0} \frac{\|g(A+B) - g(A) - \sum_{j=1}^k \frac1{j!}F_A^{(j)} (B , \dots, B)\|_\infty}{\|B\|_\infty}=0.$$
The multilinear operator $F_A^{(k)}$ is called the $k^{\mathrm{th}}$ Fr\'echet derivative of $g$ at $A$.
\end{dfn}

The following lemma describes the relationship between the derivatives introduced in Definitions~\ref{G} and~\ref{F}. 

\begin{lem}\label{GF}
If $f : \Rl \to \Rl$ is a Schwartz function and $k \in \N$, then:

(i) for every $A, B\in \mathcal {SA}$ the derivative $\frac {d^k}{dt^k} \left[ f\left( A +tB \right) \right]|_{t = 0}$ exists;

(ii) the $k^{\mathrm{th}}$ Fr\'echet derivative of
the function $A \mapsto  f(A)$ from $\mathcal {SA}$ to $\mathcal {SA}$ exists;

(iii) for every $A, B\in \mathcal {SA}$ we have 
$$F_A^{(k)} (\underbrace{B , \dots, B}_{k-\text{times}})=\frac {d^k}{dt^k} \left[ f\left( A +tB \right) \right]|_{t = 0}.$$
\end{lem}
\begin{proof}
 By~\cite[Theorem 5.7]{ACDS} for every Schwartz function $f$ the function
 $A \mapsto  f(A)$ has Fr\'echet derivatives of any order, which proves (ii).

It is shown in~\cite[Chapter VIII, \S 6 (3)]{LuSo} that the derivative \\ $\frac {d^k}{dt^k} \left[ f\left( A +tB \right) \right]|_{t = 0}$ exists and
 $$F_A^{(k)} (\underbrace{B , \dots, B}_{k-\text{times}})=\frac {d^k}{dt^k} \left[ f\left( A +tB \right) \right]|_{t = 0},$$
 which proves the first and the third assertions.
\end{proof}

%

The following theorem is proved in~\cite{PSS} (see Theorem 1.1).

\begin{thm}
\label{SSF_PSS}
Let $n\in\N$. Let $H, V$ be self-adjoint operators such that $V \in \mL_n$.
For every Schwartz function $f:\Rl \to \Rl$ the operator
$$
   f\left( H+V \right) - \sum_{k =
      0}^{n-1} \frac 1{k!} \, \frac {d^k}{dt^k} \left[ f\left( H +tV \right) \right] \biggr|_{t = 0}
$$
belongs to $\mL_1$ and
there is a unique function~$\eta_n \in L_1(\Rl)$ depending only on $n,H,V$ such that
$$
      {\rm Tr} \left( f\left( H+V \right) - \sum_{k =
      0}^{n-1} \frac 1{k!} \, \frac {d^k}{dt^k} \left[ f\left( H +tV \right) \right] \biggr|_{t = 0} \right) = \int_{\Rl}f^{(n)} (t)\, \eta_n (t)\, dt.
$$
\end{thm}

Theorem~\ref{M2} extends Theorem~\ref{SSF_PSS} to the case of weak $\mL_p$-spaces.
It is one of the main results of this paper.
The main technical ingredient of the proof of Theorem~\ref{M2} is the following estimate.

\begin{thm}\label{MainTheoremCor}
  Let $n\in\N$ and let $H, V\in B(\sH)$ be self-adjoint operators
  such that $V \in \mL_{n,\infty}$.
  For every Schwartz function $f$ we have that:

  1. The operator
  $\frac {d^n}{dt^n} \left[ f\left( H +tV\right) \right] \biggr|_{t = 0}$
  belongs to $\mL_{1,\infty}$;

2. There is a constant $c_n$ depending only on $n$ and $H$ such that for every bounded trace $\tau$ on $\mL_{1,\infty}$ the following estimate
\begin{equation}
\label{TextPrincipalEstimate}
\left|\tau \left(\frac {d^n}{dt^n} \left[ f\left( H +tV\right) \right] \biggr|_{t = 0} \right) \right| \leq \, c_n\, \left\|
        f^{(n)} \right\|_{L_\infty} \, \left\| V\right\|_{\mL_{n,\infty}}^n
\end{equation} holds.
\end{thm}

Before we proceed with the proof of Theorem~\ref{MainTheoremCor} we show that
the estimate~\eqref{TextPrincipalEstimate} easily proves Theorem~\ref{M2}.

\begin{proof}[Proof of Theorem~\ref{M2}]
Since $V \in \mL_{n,\infty} \subset \mL_{n+1}$, it follows from Theorem~\ref{SSF_PSS} that for every Schwartz function $f$, we have
$$f\left( H+V \right) - \sum_{k = 0}^{n} \frac 1{k!} \, \frac {d^k}{dt^k} \left[ f\left( H+tV\right) \right] \biggr|_{t = 0} \in \mL_1.$$

Since every bounded trace $\tau$ on $\mL_{1,\infty}$ vanishes on $\mL_1$, it follows that
$$\tau \left(f\left( H+V \right) - \sum_{k = 0}^{n-1} \frac 1{k!} \, \frac {d^k}{dt^k} \left[ f\left( H+tV\right) \right] \biggr|_{t = 0} \right)
= \frac 1{n!} \, \tau \left(\frac{d^n}{dt^n}\left[f\left(H+tV\right)\right]\biggr|_{t = 0} \right),$$
for every bounded trace $\tau$ on $\mL_{1,\infty}$.

Hence, the inequality~(\ref{TextPrincipalEstimate}) yields
$$\left|\tau \left(f\left( H+V \right) - \sum_{k = 0}^{n-1} \frac 1{k!} \, \frac {d^k}{dt^k} \left[ f\left( H+tV\right) \right] \biggr|_{t = 0} \right) \right| \leq \, c_n\, \left\|
        f^{(n)} \right\|_{L_\infty} \, \left\| V\right\|_{\mL_{n,\infty}}^n$$

Therefore, for every Schwartz functions $f$ the functional
  $$f^{(n)} \mapsto \tau\left(f\left( H+V \right) - \sum_{k = 0}^{n-1} \frac 1{k!} \, \frac {d^k}{dt^k} \left[ f\left( H+tV\right) \right] \biggr|_{t = 0}\right)$$
  is bounded in the norm of $C_b^{n}$. 
Let us consider the subspace $E \subset C(\Rl)$ consisting of all bounded functions $h$ such that $h=f^{(n)}$ for some Schwartz function $f$. 
Observe that such a function $f$ is necessarily unique, and therefore the mapping 
$$h\mapsto \varphi(h):=\tau\left(f\left( H+V \right) - \sum_{k = 0}^{n-1} \frac 1{k!} \, \frac {d^k}{dt^k} \left[ f\left( H+tV\right) \right] \biggr|_{t = 0}\right),\quad h\in E$$
is a correctly defined linear functional on $E$ satisfying 
$$|\varphi(h)|\leq  c'_n\|h\|_{C(\Rl)},\quad h\in E.$$
Hence, applying the Hahn-Banach theorem one can extend the functional $\varphi$ to a bounded functional on the space $C(\Rl)$.

By~\cite[Theorem IV.14, p.108]{RS_v1} every linear functional on the space $C(\Rl)$ of continuous functions is a linear combination of positive linear functionals.
On the other hand, by the Riesz representation theorem~\cite[436K]{Fremlin_MT_v4} for every positive linear functional $g$ on $C(\Rl)$
there exists a unique totally finite Radon measure $\eta$ on $\Rl$ such that
$$g(u) =\int_\Rl u \ d\eta, \ \forall  u \in C(\Rl).$$

Combining these two results, for every bounded trace $\tau$ on $\mL_{1,\infty}$
we deduce the existence of a unique totally finite Radon signed measure $m_n$ on $\Rl$
%
such that
$$
    \tau \left( f\left( H+V \right) - \sum_{k = 0}^{n-1} \frac 1{k!} \, \frac {d^k}{dt^k} \left[ f\left( H+tV\right) \right] \biggr|_{t = 0} \right) = \int_{\Rl} f^{(n)}(t)\, dm_n (t).
$$

Moreover, it follows from~\eqref{TextPrincipalEstimate} and the latter formula that the total variation of $m_n$ is bounded as follows: $\|m_n\| \le \|V\|^n_{\mL_{n,\infty}}$.
\end{proof}

Theorem~\ref{M1} is proved in a similar fashion, using the estimate which is analogous to that in Theorem~\ref{MainTheoremCor} (see Theorem~\ref{MR_DM} below).

The rest of the paper is devoted to the proof of Theorem~\ref{MainTheoremCor}.
To prove it we need the concept of multiple operator integrals (MOI) and the method developed in~\cite{PSS} (see proof of Theorem 5.3 there).
Two obstacles arise here. First is the construction of MOI on the quasi-Banach spaces $\mL_{p,\infty}$, $p\ge1$.
Note that in~\cite{ACDS} MOI have been constructed on the Banach spaces $\mL_p$, $p\ge1$ (using the Bochner integral) and one can interpolate them to $\mL_{p,\infty}$, $p\ge1$.
Here the second problem comes into play.
The crucial technical feature of the method from~\cite{PSS} is that at some stage one needs to swap a trace and an integral.
Which is not a problem if we are in the setting of $\mL_1$ and the standard normal trace.
(The proof of this fact, which can be found in~\cite[Theorem 3.10]{ACDS}, significantly relies on the normality of the trace.)

To overcome these obstacles in Section~\ref{sec:MOI} we construct a (Riemann) integral on the quasi-Banach spaces $\mL_{p}$, $p<1$
as a limit of integral sums of a specific form (see~\eqref{S_m} below)
and then interpolate it to $\mL_{p,\infty}$, $p\ge1$.
In this construction we follow the ideas of Turpin and Waelbroeck~\cite{TW} further developed by Kalton in~\cite{K86b}.
Using these new integrals we introduce the notion of MOI in $\mL_{p,\infty}$, $p\ge1$
and prove that they are well-defined for a sufficiently large class of functions.
Also this construction of the integral allows us to swap a trace (not necessarily normal) and an integral for sufficiently large class of integrands (see Theorem~\ref{swap} below).

Next, in sections~\ref{sec:proof} and~\ref{sec:DM} we employ multiple operator integrals to prove Theorems~\ref{M2} and~\ref{M1} respectively.
In Section~\ref{sec:Unbounded} we also use the multiple operator integrals to establish the perturbation formulae in the case of an unbounded operator $H$.
The last section is an appendix, where we have gathered the results concerning the polylinear interpolation which are used in this paper.

The following Riesz-Fischer type theorem is used in Theorem~\ref{int} below.
Although this result is folklore, we give a short proof for the convenience of a reader.

\begin{thm}\label{RF} Let $0<p< 1$.
 If $A_k \in \mL_p$ are such that $\sum_{k=0}^\infty \|A_k\|_{\mL_p}^p <\infty$,
 then $\sum_{k=0}^\infty A_k$ converges in $\mL_p$ and
 $$\|\sum_{k=0}^\infty A_k\|_{\mL_p}^p \le \sum_{k=0}^\infty \|A_k\|_{\mL_p}^p.$$
\end{thm}

\begin{proof}
 Since the quasi-norm of $\mL_p$ satisfies the following inequality (see e.g.~\cite[Theorem 4.9]{FK})
\begin{equation}\label{qin}
 \|A+B\|_{\mL_p}^p \le \|A\|_{\mL_p}^p +\|B\|_{\mL_p}^p, \ A,B \in \mL_p,
\end{equation}
it follows that
$$\|\sum_{k=m}^n A_k\|_{\mL_p}^p \le \sum_{k=m}^n \|A_k\|_{\mL_p}^p \to 0, \ \text{as} \ m,n \to \infty.$$

Hence, $\{ \sum_{k=0}^n A_k \}_{n=1}^\infty$ is a Cauchy sequence in $\mL_p$ with respect to the metric $d_p$ defined as follows
\begin{equation}\label{d_p}
 d_p(A,B) = \|A-B\|_{\mL_p}^p.
\end{equation}

Since $\mL_p$ is complete with respect to the metric $d_p$ (see e.g.~\cite{S_indag}), it follows that $\sum_{k=0}^\infty A_k$ converges in $\mL_p$.

Next, using the well-known Fatou property of $\mL_p$, we obtain
$$\|\sum_{k=0}^\infty A_k\|_{\mL_p}^p = \lim_{n\to \infty}\|\sum_{k=0}^n A_k\|_{\mL_p}^p \le \lim_{n\to \infty}\sum_{k=0}^n \|A_k\|_{\mL_p}^p = \sum_{k=0}^\infty \|A_k\|_{\mL_p}^p.$$
\end{proof}

 Note that the constants below may vary from line to line and even within the line, although the same letter $C$ is used to denote them.
We do this when the value of constants has no relevance
to statements or proofs.

\section{Multiple operator integrals}\label{sec:MOI}
We fix an integer $n\ge 1$ and denote by $\Rl_+^{n+1}$ the positive cone of $\Rl^{n+1}$.
Let~$\cC_n$ be the class of functions~$\phi: \Rl^{n + 1} \mapsto \Cx$
admitting the representation
\begin{equation}\label{Crep}
  \phi (\lambda_0, \ldots, \lambda_n) =
  \int_{\Rl_+^{n+1}} \prod_{j = 0}^n a_j (\lambda_j, s)\, d\nu(s),
\end{equation}
for some
bounded
Borel functions $ a_{j} \left( \cdot, s \right) : \Rl \mapsto
\Cx $ and measure $\nu$ on $\Rl_+^{n+1}$ such that $|\nu|$ is finite and
\begin{equation}\label{phi_norm}
 \int_{\Rl_+^{n+1}} \prod_{j = 0}^n \left\| a_j(\cdot, s) \right\|_\infty \,  d| \nu |(s)<\infty.
\end{equation}

Following~\cite[Definition 4.1]{ACDS} for every function $\phi\in\cC_n$ we define
the corresponding operator integral.

\begin{dfn}
  \label{MOI}
 Let $1\leq p_j \leq \infty$, $1\leq j\leq n$ and $p$ be such that $\frac 1{p} = \sum_{j=1}^n\frac 1{p_j}$.
 For every~$\phi \in \cC_n$ we define the operator
 $$T_\phi : \mL_{p_1} \times \dots \times \mL_{p_n} \to \mL_{p}$$
 as follows:
  \begin{equation}
    \label{ClassNewDefRep}
    \begin{aligned}
     T_\phi (x_1, x_2, \ldots, x_n)
    := &\int_{\Rl_+^{n+1}} a_0(H, s)\, x_1 \, a_1 (H, s)\, x_2 \cdot \ldots\cdot x_n a_n (H, s)\, d \nu(s),
    \end{aligned}
  \end{equation}where the integral above is the Bochner integral on $\mL_{p}$.
\end{dfn}

For the rest of the section we prove auxiliary results required to extend the notion of
multiple operator integrals to the case of $\mL_{p,\infty}$ spaces.

First we define the integration in an arbitrary quasi-Banach space $(X, \| \cdot \|_X)$.
However, we will only use it in the case when $X$ either $\mL_{p,\infty}$ or $\mL_{p}$.

For every $k_0, k_1, \dots, k_n \in \Rl^{n+1}$ we denote $\z{k}:=(k_0, k_1, \dots, k_n)$, $\frac{\z k}{2^{m}} := ( \frac{k_0}{2^{m}}, \frac{k_1}{2^{m}}, \dots, \frac{k_n}{2^{m}})$
and $\frac{2\lfloor \frac{\z k}2 \rfloor}{2^{m}} := (\frac{2\lfloor \frac{k_0}2 \rfloor}{2^{m}}, \frac{2\lfloor \frac{k_1}2 \rfloor}{2^{m}},\dots, \frac{2\lfloor \frac{k_n}2 \rfloor}{2^{m}}).$
Also, by writing $\z k\ge0$ we understand, that all $k_0, k_1, \dots, k_n \ge0$.
\begin{dfn}\label{op_int}
Let $(X, \| \cdot \|_X)$ be a quasi-Banach space.
Let $n\in\N$.
Let
$u : \Rl_{+}^{n+1} \to X$.
The integral sums $S_m$, $m\ge1$ are defined by the following formula
\begin{equation}\label{S_m}
 S_m := \frac1{2^{m(n+1)}} \sum u( \frac{\z k}{2^{m}}),
\end{equation}
where the sum is taken over all $0\le \z k \in \mathbb Z^{n+1}$.

If the series in~\eqref{S_m} and the sequence $\{S_m\} \subset X$ itself are convergent in the quasi-norm $\| \cdot \|_X$, then we set
$$\int_{\Rl_{+}^{n+1}} u(s) \ ds := \lim_{m\to \infty} S_m.$$
\end{dfn}

Throughout the paper all integrals are understood in the sense of the preceding definition, unless explicitly specified.

Note that, in particular, the latter definition introduces a notion of an integral of $\mL_p$-valued functions, $p\ge1$.
However, for $p\ge1$ the space $\mL_p$ is Banach. Therefore, in these settings, it makes sense also to speak of the Bochner integral of an $\mL_p$-valued function.
So, for $p\ge1$ there are two notions of an integral of an $\mL_p$-valued function: the Bochner integral and the integral in the sense of Definition~\ref{op_int}.
Theorem~\ref{bohner} below shows that if $p\ge1$ and a function $u: \Rl_{+}^{n+1} \to \mL_p$ (satisfying some additional conditions) is integrable in the sense of Bochner, then
$u$ is also integrable in the sense of Definition~\ref{op_int} and the Bochner integral coincides with that introduced in Definition~\ref{op_int}.

We start with a result proving that the integral (in the sense of Definition~\ref{op_int}) exists for a wide class of functions $u : \Rl_{+}^{n+1} \to \mL_p$.

\begin{thm}\label{int} Let $n\in\N$ and $p> \frac{n+1}{n+2}$.
Let
$u : \Rl^{n+1} \to \mL_p$ be such that $u=wv$, where
$w : \Rl^{n+1} \to \mL_p$ is a bounded Lipschitz function and
$v : \Rl^{n+1} \to \Cx$ is such that
\begin{equation}\label{sw1}
 |v(s)|, |\triangledown v(s)| \le  \frac{C_\alpha}{(1+\sum_{i=0}^n |s_i|)^{\alpha}}
\end{equation}
for every $s=(s_0, \dots, s_n) \in \Rl^{n+1}$ and every $\alpha >0$,
where $\triangledown$ denotes the gradient of a function and $C_\alpha$ is a constant depending on $\alpha$ only.

We have

1. For every $m \in \N$ the series
$$S_m:=\frac1{2^{m(n+1)}} \sum_{\z{k} \ge 0} u( \frac{\z k}{2^{m}})$$
is convergent in $\mL_p$ and $\|S_m\|_{\mL_p} \le C \cdot 2^{m(n+1)/p}$;

2. The sequence $\{S_m\}$ is convergent in $\mL_p$ and, so, the integral $\int_{\Rl_+^{n+1}} u(s) \ ds$ is defined
and $\|\int_{\Rl_+^{n+1}} u(s) \ ds\|_{\mL_p} \le C \max\{\|w\|_{Lip}, \|w\|_{L_\infty(\mL_p)}\}$;

3. The following estimate holds: 
$$\| S_m - \int_{\Rl_+^{n+1}} u(s) \ ds \|_{\mL_p} \le C \max\{\|w\|_{Lip}, \|w\|_{L_\infty(\mL_p)}\} 2^{m(\frac{n+1}p-(n+2))}, \ \forall m\ge1,$$
where $\|w\|_{L_\infty(\mL_p)}:= \sup_{s \in \Rl^{n+1}} \|w(s)\|_{\mL_p} <\infty$.
\end{thm}

\begin{proof}
We prove the theorem for the case $\frac{n+1}{n+2}<p<1$.
For $p\ge1$ the proof is similar (and easier), with the only difference that instead of Theorem~\ref{RF}
we use the triangle inequality for $\mL_p$-norm.

1. For every $m\in \N$ we first prove that the series
$$\sum_{\z k \ge0} \|u(\frac{\z k}{2^{m}})\|_{\mL_p}^p$$
is convergent.

Using the first inequality~\eqref{sw1} for $v$ with $\alpha=\frac{n+2}p$ and the fact that the function $w$ is bounded, we obtain 
\begin{align*}
\sum_{\z k \ge0} \|u(\frac{\z k}{2^{m}})\|_{\mL_p}^p
 &\le \sum_{\z k \ge0} \|w\|_{L_\infty(\mL_p)}^p \cdot |v|^p(\frac{\z k}{2^{m}})
 \le C_\alpha  \sum_{\z k \ge0} \frac{1}{(1+2^{-m}\sum_{i=0}^n k_i)^{n+2}}.
\end{align*}

For the latter sum we have
\begin{equation}\label{sum}
\begin{aligned}
 \sum_{\z k \ge0} \frac{1}{(1+2^{-m}\sum_{i=0}^n k_i)^{n+2}} &= \sum_{l =0}^\infty \frac{1}{(1+\frac{l}{2^{m}})^{n+2}} \sum_{k_0+\dots + k_n =l} 1
 \le \sum_{l =0}^\infty \frac{l^n}{(1+\frac{l}{2^{m}})^{n+2}}\\
 &= 2^{mn} \sum_{l =0}^\infty \frac{\left( \frac{l}{2^{m}} \right)^n}{(1+\frac{l}{2^{m}})^{n+2}}
 \le 2^{mn} \sum_{l =0}^\infty \frac{\left( \lfloor \frac{l}{2^{m}} \rfloor +1\right)^n}{(1+\lfloor \frac{l}{2^{m}} \rfloor)^{n+2}}\\
 &= 2^{m(n+1)} \sum_{j=0}^\infty \frac{(j +1)^n}{(1+j)^{n+2}}
 = 2^{m(n+1)} \frac{\pi^2}6.
\end{aligned}
\end{equation}

Hence, for every $m \in \N$ Theorem~\ref{RF} yields
 $$ \|S_m\|_{\mL_p}^p \le \frac1{2^{mp(n+1)}}\sum_{\z k \ge0} \|u(\frac{\z k}{2^{m}})\|_{\mL_p}^p \le 2^{m(1-p)(n+1)} \frac{\pi^2}6 <\infty,$$
 that is series $S_m$ is convergent in $\mL_p$.

2. For every $m \ge2$ consider
$$S_{m-1} = \frac1{2^{(m-1)(n+1)}} \sum_{\z k\ge0} u(\frac{\z k}{2^{m-1}}).$$
By the first part of the proof the latter series is absolutely convergent in $\mL_p$, so one can change the order of summation.
We have
$$S_{m-1} = \frac1{2^{(m-1)(n+1)}} \sum_{\z k\ge0} u(\frac{2 \z k}{2^{m}})
= \frac1{2^{m(n+1)}} \sum_{\z k\ge0} u(\frac{2\lfloor \frac{\z k}2 \rfloor}{2^{m}}),$$
since in the latter sum every element is repeated $2^{n+1}$-times.

Hence, for $m\ge 2$ we have
\begin{equation}\label{es10}
\begin{aligned}
 S_m-S_{m-1} &= \frac1{2^{m(n+1)}} \sum_{\z k\ge0} \left[ u(\frac{\z k}{2^{m}}) - u(\frac{2\lfloor \frac{\z k}2 \rfloor}{2^{m}}) \right]\\
 &= \frac1{2^{m(n+1)}} \sum_{\z k\ge0} \left[ w(\frac{\z k}{2^{m}})v(\frac{\z k}{2^{m}})
- w(\frac{2\lfloor \frac{\z k}2 \rfloor}{2^{m}}) v(\frac{2\lfloor \frac{\z k}2 \rfloor}{2^{m}}) \right]\\
 &= \frac1{2^{m(n+1)}} \sum_{\z k\ge0} \left( w(\frac{\z k}{2^{m}})- w(\frac{2\lfloor \frac{\z k}2 \rfloor}{2^{m}})  \right) v(\frac{\z k}{2^{m}})\\
 &+ \frac1{2^{m(n+1)}} \sum_{\z k\ge0} w(\frac{2\lfloor \frac{\z k}2 \rfloor}{2^{m}})
 \left( v(\frac{\z k}{2^{m}})- v(\frac{2\lfloor \frac{\z k}2 \rfloor}{2^{m}}) \right).
\end{aligned}
\end{equation}

For the first sum, since $w$ is Lipschitz and $v$ satisfies~\eqref{sw1} (with $\alpha=\frac{n+2}p$), we have
\begin{equation}\label{es11}
\begin{aligned}
&\sum_{\z k\ge0} \left\|\left( w(\frac{\z k}{2^{m}})- w(\frac{2\lfloor \frac{\z k}2 \rfloor}{2^{m}})  \right) v(\frac{\z k}{2^{m}}) \right\|_{\mL_p}^p\\
&\le \sum_{\z k\ge0}  \left\|w(\frac{\z k}{2^{m}})- w(\frac{2\lfloor \frac{\z k}2 \rfloor}{2^{m}})\right\|_{\mL_p}^p \cdot |v|^p(\frac{\z k}{2^{m}})\\
&\le \|w\|_{Lip}^p \ 2^{-mp} \sum_{\z k\ge0}  \frac{C_\alpha}{(1+2^{-m}\sum_{i=0}^n |k_i|)^{n+2}} \\
&\le  C_\alpha \ \|w\|_{Lip}^p \ 2^{m(n+1-p)} ,
\end{aligned}
\end{equation}
due to~\eqref{sum}.

For the second sum, since $w$ is bounded on $\Rl^{n+1}$ and $\triangledown v$ satisfies~\eqref{sw1} (with $\alpha=\frac{n+2}p$), we have
\begin{equation}\label{es12}
\begin{aligned}
&\sum_{\z k\ge0} \left\|w(\frac{2\lfloor \frac{\z k}2 \rfloor}{2^{m}})
 \left( v(\frac{\z k}{2^{m}})- v(\frac{2\lfloor \frac{\z k}2 \rfloor}{2^{m}}) \right)\right\|_{\mL_p}^p\\
&\le  \|w\|_{L_\infty(\mL_p)}^p \sum_{\z k\ge0} \left|v(\frac{\z k}{2^{m}})- v(\frac{2\lfloor \frac{\z k}2 \rfloor}{2^{m}})\right|^p\\
&\le  \|w\|_{L_\infty(\mL_p)}^p \sum_{\z k\ge0} \frac{C_\alpha 2^{-mp}}{(1+2^{-m}\sum_{i=0}^n |k_i|)^{n+2}}\\
&\le C_\alpha \|w\|_{L_\infty(\mL_p)}^p \ 2^{m(n+1-p)},
\end{aligned}
\end{equation}
due to~\eqref{sum}.

Combining~\eqref{es10} with~\eqref{es11} and~\eqref{es12} and using Theorem~\ref{RF} we obtain
\begin{equation}\label{es13}
\begin{aligned}
 \|S_m-S_{m-1}\|_{\mL_p}^p
 &\le \frac1{2^{mp(n+1)}} \sum_{\z k\ge0} \left\|\left( w(\frac{\z k}{2^{m}})- w(\frac{2\lfloor \frac{\z k}2 \rfloor}{2^{m}})  \right) v(\frac{\z k}{2^{m}})\right\|_{\mL_p}^p \\
 &+ \frac1{2^{mp(n+1)}} \sum_{\z k\ge0} \left\|w(\frac{2\lfloor \frac{\z k}2 \rfloor}{2^{m}}) \left( v(\frac{\z k}{2^{m}})- v(\frac{2\lfloor \frac{\z k}2 \rfloor}{2^{m}}) \right)\right\|_{\mL_p}^p\\
 &\le \frac{C \ 2^{m(n+1-p)}}{2^{mp(n+1)}}(\|w\|_{Lip}^p + \|w\|_{L_\infty(\mL_p)}^p) \\
 &\le C \max\{\|w\|_{Lip}, \|w\|_{L_\infty(\mL_p)}\}^p \ 2^{m(n+1-(n+2)p)} .
\end{aligned}
\end{equation}

Since $p > \frac{n+1}{n+2}$, it follows that $n+1-(n+2)p<0$. Therefore,
$$\sum_{m=2}^\infty \|S_m-S_{m-1}\|_{\mL_p}^p < \infty.$$

Hence, $\{ S_m \}_{m=1}^\infty$ is a Cauchy sequence in $\mL_p$ with respect to the metric $d_p$ defined by the formula~\eqref{d_p}.
Since $\mL_p$ is complete with respect to the metric $d_p$, it follows that
the limit of the sequence $\{ S_m \}_{m=1}^\infty$ exists and, so, the integral $\int_{\Rl_+^{n+1}} u(s) \ ds$ is well-defined in the sense of Definition~\ref{op_int}.

3. The inequality~\eqref{es13} yields
\begin{align*}
 \| S_m - \int_{\Rl_+^{n+1}} u(s) \ ds \|_{\mL_p}^p &\le \sum_{k=m}^\infty \|S_k-S_{k-1}\|_{\mL_p}^p \\
 &\le C \max\{\|w\|_{Lip}, \|w\|_{L_\infty(\mL_p)}\}^p \sum_{k=m}^\infty 2^{k(n+1-(n+2)p)}.
\end{align*}

Hence,
$$\| S_m - \int_{\Rl_+^{n+1}} u(s) \ ds \|_{\mL_p}  \le C \max\{\|w\|_{Lip}, \|w\|_{L_\infty(\mL_p)}\} 2^{m(\frac{n+1}p-(n+2))}.$$

Note that for $p\ge1$ we can prove a stronger estimate.
Indeed, using~\eqref{sw1} with $\alpha=n+2$ we obtain the following versions of~\eqref{es11} and~\eqref{es12}:
\begin{equation}\label{em0}
\sum_{\z k\ge0} \left\|\left( w(\frac{\z k}{2^{m}})- w(\frac{2\lfloor \frac{\z k}2 \rfloor}{2^{m}})  \right) v(\frac{\z k}{2^{m}}) \right\|_{\mL_p}
\le  C_\alpha \ \|w\|_{Lip} \ 2^{mn}
\end{equation}
and
\begin{equation}\label{em1}
\sum_{\z k\ge0} \left\|w(\frac{2\lfloor \frac{\z k}2 \rfloor}{2^{m}})
 \left( v(\frac{\z k}{2^{m}})- v(\frac{2\lfloor \frac{\z k}2 \rfloor}{2^{m}}) \right)\right\|_{\mL_p}
\le C_\alpha \|w\|_{L_\infty(\mL_p)} \ 2^{mn}.
\end{equation}

Hence, combining~\eqref{es10} with~\eqref{em0}, \eqref{em1}
and using the triangle inequality for the norm (instead of the inequality~\eqref{qin} for the quasi-norm) we obtain
$$\|S_m-S_{m-1}\|_{\mL_p}\le C \max\{\|w\|_{Lip}, \|w\|_{L_\infty(\mL_p)}\}^p \ 2^{-m}$$
and
$$\| S_m - \int_{\Rl_+^{n+1}} u(s) \ ds \|_{\mL_p}  \le C \max\{\|w\|_{Lip}, \|w\|_{L_\infty(\mL_p)}\} 2^{-m}.$$

\end{proof}

Now we show that our Definition~\ref{op_int} coincides with the definition of the Bochner integral in $\mL_p$, $p\ge1$.
We need the following result established in~\cite{vNeerven}.

Recall that a measure $\nu$ on a topological space $\Omega$ is called tight if for every $\varepsilon>0$
there is a compact set $K \subset \Omega$ such that $|\nu | (\Omega \setminus K) < \varepsilon$.

\begin{thm}[\cite{vNeerven}]\label{vN}
 Let $\nu$ be a tight Borel measure on a metric space $\Omega$, let $X$ be
a Banach space. If $w : (\Omega, \nu) \to X$ is Bochner integrable,
then for every sequence of partitions
$$P^{(m)} = \{ \Omega_j^{(m)}\}_{j\ge0}$$
of $\Omega$ satisfying $\lim_{m\to \infty} \sup_{j\ge0} {\rm diam} \ \Omega_j^{(m)} = 0$
there exists a sequence of sample point sets
$$Z^{(m)} = \{ z_j^{(m)} \in \Omega_j^{(m)} \}_{j\ge0}$$
such that
$$\lim_{m\to \infty} \left\| \sum_{j=0}^\infty w(z_j^{(m)}) \nu (\Omega_j^{(m)}) - \int_{\Omega} w(s) \ d\nu(s) \right\|_X= 0,$$
where the latter integral is the Bochner integral in $X$.
\end{thm}

The following result shows that our new definition of integral extends the notion of Bochner integral to $\mL_{p,\infty}$.
\begin{thm}\label{bohner}
Let $n\in\N$ and $p\ge1$.
Let
$u : \Rl^{n+1} \to \mL_p$ be such that $u=wv$, where
$w : \Rl^{n+1} \to \mL_p$ is a bounded Lipschitz function and
$v : \Rl^{n+1} \to \Cx$ satisfies~\eqref{sw1}.
Then $\lim_{m\to\infty} S_m$ ($S_m$ is defined as in~\eqref{S_m}) is equal to the Bochner integral $\int_{\Rl_+^{n+1}} w(s) \ d\nu(s)$, where $d\nu(s)= v(s) ds$.
\end{thm}

\begin{proof}
In this proof all integrals are understood in the sense of Bochner.

Set
$$P^{(m)} = \{ [\frac{k_0}{2^{m}}, \frac{k_0+1}{2^{m}}] \times \dots \times [\frac{k_n}{2^{m}}, \frac{k_n+1}{2^{m}}],
\quad \forall \ \z k \ge0\}.$$
Note that the diameter of every set from this partition is $\sqrt{\frac{n+1}{2^m}}\to 0$, $m\to\infty$.
Also note, that since $v$ satisfies~\eqref{sw1} the measure $\nu$ is a finite tight Borel measure on $\Rl^{n+1}$
and the function $w$ is Bochner integrable with respect to the measure $\nu$ (since $w$ is bounded Lipschitz and $\nu$ is finite).

Hence, by Theorem~\ref{vN} for every $\z k \ge0$ there exists $z^{(m)}(\z k) \in [\frac{k_0}{2^{m}}, \frac{k_0+1}{2^{m}}] \times \dots \times [\frac{k_n}{2^{m}}, \frac{k_n+1}{2^{m}}],$ such that
$$\left\| \frac1{2^{m(n+1)}} \sum_{\z k\ge0} w(z^{(m)}(\z k)) v(z^{(m)}(\z k))- \int_{\Rl_+^{n+1}} w(s) \ d\nu(s) \right\|_{\mL_p} \to 0, m\to\infty.$$

Similarly to~\eqref{em0} and~\eqref{em1} (using~\eqref{sw1} with $\alpha=n+2$) we obtain
\begin{equation}\label{ooo}
\sum_{\z k\ge0} \left\|w(\frac{\z k}{2^{m}})v(\frac{\z k}{2^{m}}) - w(z^{(m)}(\z k))v(z^{(m)}(\z k))\right\|_{\mL_p} \le C 2^{mn}.
\end{equation}
Here $C$ is a constant depending on $w$ and $n$.

Next, since $\mL_p$, $p\ge1$ is a Banach space, it follows that
\begin{align*}
&\left\| \frac1{2^{m(n+1)}} \sum_{\z k\ge0} w(\frac{\z k}{2^{m}})v(\frac{\z k}{2^{m}}) - \int_{\Rl_+^{n+1}} w(s) \ d\nu(s) \right\|_{\mL_p}\\
\le
&\left\| \frac1{2^{m(n+1)}} \sum_{\z k\ge0} w(z^{(m)}(\z k))v(z^{(m)}(\z k)) - \int_{\Rl_+^{n+1}} w(s) \ d\nu(s) \right\|_{\mL_p} \\
+&\left\|\frac1{2^{m(n+1)}} \sum_{\z k\ge0} [w(\frac{\z k}{2^{m}})v(\frac{\z k}{2^{m}})) - w(z^{(m)}(\z k))v(z^{(m)}(\z k))]\right\|_{\mL_p}.
\end{align*}

For the second term above, using~\eqref{ooo} and the triangle inequality we have
\begin{align*}
&\left\|\frac1{2^{m(n+1)}} \sum_{\z k\ge0} [w(\frac{\z k}{2^{m}})v(\frac{\z k}{2^{m}})) - w(z^{(m)}(\z k))v(z^{(m)}(\z k))]\right\|_{\mL_p} \\
\le &\frac1{2^{m(n+1)}} \sum_{\z k\ge0} \left\|w(\frac{\z k}{2^{m}})v(\frac{\z k}{2^{m}})) - w(z^{(m)}(\z k))v(z^{(m)}(\z k))\right\|_{\mL_p}\le  C  2^{-m}.
\end{align*}

Consequently,
$$\left\| \frac1{2^{m(n+1)}} \sum_{\z k\ge0} u(\frac{\z k}{2^{m}}) - \int_{\Rl_+^{n+1}} w(s) \ d\nu(s) \right\|_{\mL_p} \to 0, m\to\infty.$$

In other words, the Bochner integral
$$\int_{\Rl_+^{n+1}} w(s) \ d\nu(s) = \lim_{m\to\infty} S_m.$$
\end{proof}

Now we extend the notion of the multiple operator integral to the quasi-Banach ideal $\mL_{p,\infty}$.

\begin{dfn}\label{quasiMOI}
 Let $1\leq p_j \leq \infty$, $1\leq j\leq n$ and let $p$ be such that $\frac 1{p} = \sum_{j=1}^n\frac 1{p_j}$.
 Let $\phi\in \cC_n$ admits the representation~\eqref{Crep} with a measure $\nu$ being absolute continuous.
 A function $u_\phi : \Rl_+^{n+1} \to \mL_{p,\infty}$ is defined as follows:
\begin{equation}\label{u_dfn}
 u_\phi(s) := a_0(H,s)x_1 a_1(H,s)x_2 \cdots x_n a_n(H,s) \nu'(s), \ s \in \Rl_+^{n+1},
\end{equation}
where the function $\nu' : \Rl_+^{n+1} \to \Cx$ is the Radon-Nikodym derivative of $\nu$.

 If the integral $\int_{\Rl_+^{n+1}} u_\phi(s) ds \ $
 exists in the sense of Definition~\ref{op_int},
 then  we define the multilinear operator
 $$T_{\phi} : \mL_{p_1,\infty} \times \dots \times \mL_{p_n,\infty} \to \mL_{p,\infty}$$
by setting
\begin{equation}\label{new_int}
 T_{\phi}(x_1, \dots, x_n) := \int_{\Rl_+^{n+1}} u_\phi(s) \ ds.
\end{equation}
\end{dfn}

Theorem~\ref{qn_conv} below describes the class of functions $\phi\in \cC_n$ for which the latter definition makes sense.
We first prove the following auxiliary result.

\begin{lem}\label{qn_estimate}
Let $1\leq p_j \leq \infty$, $1\leq j\leq n$ and let $p$ be such that $\frac 1{p} = \sum_{j=1}^n\frac 1{p_j}$.
Let $x_j \in \mL_{p_j}$ and let~$H\in B(\sH)$.
 If the functions $a_j(H,\cdot): \Rl^{n+1} \to B(\sH)$ are bounded and Lipschitz for all $0\le j \le n$,
 then the function
 $w : \Rl^{n+1} \to \mL_p$ given by
 $$w(s) := a_0(H,s)x_1 a_1(H,s)x_2 \cdots x_n a_n(H,s)$$
 is bounded, Lipschitz and
 $$\|w\|_{Lip}, \ \|w\|_{L_\infty(\mL_p)} \le  C \prod_{j=1}^n\|x_j\|_{\mL_{p_j}}.$$
\end{lem}

\begin{proof}
 We shall prove this lemma for $n=2$. The case of $n\ge 3$ can be treated similarly.

Using H\"older's inequality we obtain
$$\|w\|_{L_\infty(\mL_p)} \le \prod_{j=0}^{n}\|a_j\|_{L_\infty} \cdot \prod_{j=1}^{n}\|x_j\|_{\mL_{p_j}}.$$

For $s, t \in \Rl^{3}$ we have
\begin{align*}
 w(s) - w(t) &= a_0(H,s)x_1 a_1(H,s)x_2 a_2(H,s)  - a_0(H,t)x_1 a_1(H,t)x_2 a_2(H,t) \\
 &= [a_0(H,s) - a_0(H,t)] x_1 a_1(H,s)x_2 a_2(H,s) \\
 &+ a_0(H,t) x_1 [a_1(H,s) - a_1(H,t)] x_2 a_2(H,s) \\
 &+ a_0(H,t)x_1 a_1(H,t) x_2 [a_2(H,s) - a_2(H,t)].
\end{align*}

Using the quasi-triangle inequality and H\"older's inequality, we obtain
\begin{align*}
&\| w(s) - w(t)\|_{\mL_p}\\
&\le C \|x_1\|_{\mL_{p_1}} \|x_2\|_{\mL_{p_2}}
\left( \|a_0(H,s) - a_0(H,t)\|_\infty \|a_1(H,s)\|_\infty \|a_2(H,s)\|_\infty  \right.\\
 &\left. + \|a_0(H,t)\|_\infty \|a_1(H,s) - a_1(H,t)\|_\infty \|a_2(H,s)\|_\infty  \right.\\
 &\left. + \|a_0(H,t)\|_\infty \|a_1(H,t)\|_\infty \|a_2(H,s) - a_2(H,t)\|_\infty  \right).
\end{align*}

Since all functions $a_j$ are Lipschitz in $B(\sH)$, it follows that
$w$ is Lipschitz in $\mL_p$ and
$$\|w\|_{Lip} \le  C \prod_{j=1}^n\|x_j\|_{\mL_{p_j}},$$
where $C$ depends on $p$, $H$, $\|a_j\|_{Lip}$ and $\|a_j\|_\infty$, $j=0,1,2$.
\end{proof}


\begin{thm}\label{qn_conv}
 Let $1\leq p_j \leq \infty$, $1\leq j\leq n$ and let $p$ be such that $\frac 1{p} = \sum_{j=1}^n \frac 1{p_j}$.
Let $x_j \in \mL_{p_j,\infty}$ and let~$H\in B(\sH)$.
Let $\phi\in \cC_n$ be such that

1. $\phi$ admits the representation~\eqref{Crep} with a measure $\nu$ being absolute continuous
and such that its Radon-Nikodym derivative $\nu'$ satisfies~\eqref{sw1};

2. the functions $a_j(H,\cdot): \Rl^{n+1} \to B(\sH)$ from the representation~\eqref{Crep} are bounded and Lipschitz for all $0\le j \le n$.

Define a function $u_\phi : \Rl^{n+1} \to \mL_{p, \infty}$ by setting
$$u_\phi(s_0, s_1, \dots, s_n) := a_0(H,s)x_1 a_1(H,s)x_2 \cdots x_n a_n(H,s) \nu'(s), \ s \in \Rl^{n+1}.$$

 We have that
the integral $\int_{\Rl_+^{n+1}} u_\phi(s) \ ds$ is well-defined and
$$\|S_m - \int_{\Rl_+^{n+1}} u_\phi(s) \ ds \|_{\mL_{p, \infty}} \to 0,$$
where $S_m$ is defined as follows:
$$S_m=\frac1{2^{m(n+1)}} \sum_{\z k\ge 0} u_\phi( \frac{\z k}{2^{m}}).$$
\end{thm}

\begin{proof}
First note that the series in the definition of $S_m$ are convergent in the quasi-norm of $\mL_{p,\infty}$ for every $m\in \N$.
The proof of this assertion is a direct repetition of that of Theorem~\ref{int} (i).
Therefore, we omit it.

Set $$w(s) := a_0(H,s)x_1 a_1(H,s)x_2 \cdots x_n a_n(H,s), s\in \Rl^{n+1} \ \text{and} \  v=\nu'.$$
For every $\alpha>\frac{n+1}{n+2}$ and some $0 < \alpha_1, \dots, \alpha_n < \infty$ satisfying $\frac1\alpha =\frac 1 {\alpha_1} + \ldots + \frac 1 {\alpha_n}$ consider operators
$$S_m, S : \mL_{\alpha_1} \times \ldots \times \mL_{\alpha_n} \to \mL_\alpha$$
defined as follows:
$$S_m=\frac1{2^{m(n+1)}} \sum_{\z k\ge 0} u_\phi( \frac{\z k}{2^{m}}), \ m \in \N \quad \text{and} \quad S=\int_{\Rl_+^{n+1}} u_\phi(s) \ ds.$$

First note, that by Lemma~\ref{qn_estimate} the function $u_\phi$ satisfies the conditions of Theorem~\ref{int},
so by Theorem~\ref{int} the above operators are well-defined.
Also Theorem~\ref{int} yields

1. $\|S_m\|_{\mL_{\alpha_1} \times \ldots \times \mL_{\alpha_n} \to \mL_\alpha} \le C 2^{\frac{m}{\alpha} (n+1)};$

2. $\|S\|_{\mL_{\alpha_1} \times \ldots \times \mL_{\alpha_n} \to \mL_\alpha} \le C \max\{ \|w\|_{Lip}, \ \|w\|_{L_\infty(\mL_\alpha)}\};$

3. $\|S_m -S\|_{\mL_{\alpha_1} \times \ldots \times \mL_{\alpha_n} \to \mL_\alpha} \le C \max\{ \|w\|_{Lip}, \ \|w\|_{L_\infty(\mL_\alpha)}\} 2^{m(\frac{n+1}\alpha-(n+2))};$

Therefore, by Theorem~\ref{Gen_int} below operators $S_m, S, S_m-S$ are bounded from $\mL_{p_1,\infty} \times \ldots \times \mL_{p_n,\infty}$ to $\mL_{p, \infty}$ for every $m \in \N$.
In other words,  the integral $\int_{\Rl_+^{n+1}} u_\phi(s) \ ds$ is defined.
Moreover,
$$\|S_m -S\|_{\mL_{p_1, \infty} \times \ldots \times \mL_{p_n, \infty} \to \mL_p, \infty} \mathop{\longrightarrow}\limits_{m\to\infty} 0.$$

In other words,
$$\|S_m - \int_{\Rl_+^{n+1}} u_\phi(s) \ ds \|_{\mL_{p, \infty}} \to 0.$$
\end{proof}

\begin{cor}\label{int_cor}
 Let $1\leq p_j \leq \infty$, $1\leq j\leq n$ and let $p$ be such that $\frac 1{p} = \sum_{j=1}^n \frac 1{p_j}$.
  If $\phi\in \cC_n$ satisfies the assumptions of Theorem~\ref{qn_conv}, then the operator integral 
 $$T_\phi : \mL_{p_1,\infty} \times \dots \times \mL_{p_n,\infty} \to \mL_{p,\infty}$$
given in Definition~\ref{quasiMOI} is well-defined and for every $x_j \in \mL_{p_j,\infty}$, $1\leq j\leq n$ the element
 $T_\phi(x_1, \dots, x_n)$ (which is given by the integral in~\eqref{new_int}) is the limit of integral sums of the form~\eqref{S_m} with respect to the quasi-norm of $\mL_{p,\infty}$.
\end{cor}

The following result is the crucial part of the proof of Theorem~\ref{MainTheoremCor} below.

\begin{thm}\label{swap}
 Let $n\in\N$, $p=1$ and let $\phi\in \cC_n$ and $u_\phi$ be as in Theorem~\ref{qn_conv}.
For every bounded trace $\tau$ on $\mL_{1,\infty}$
we have
$$\tau (\int_{\Rl_+^{n+1}} u_\phi(s) \ ds) = \int_{\Rl_+^{n+1}} \tau(u_\phi(s)) \ ds.$$
\end{thm}

\begin{proof}
 By Theorem~\ref{qn_conv} the integral $\int_{\Rl_+^{n+1}} u_\phi(s) \ ds$ is a limit (with respect to the quasi-norm of $\mL_{1,\infty}$)
of integral sums of the form~\eqref{S_m} and since $\tau$ is continuous one can swap the trace and the integral.
\end{proof}

\begin{rem}
 The result of Theorem~\ref{swap} holds in a more general setting: instead of a trace one can take an arbitrary bounded linear functional on $\mL_{1,\infty}$.
\end{rem}

\begin{rem}
For the normal trace $\rm Tr$ on $\mL_{1}$ the proof of the latter equality is much simpler (see e.g.~\cite[Theorem 3.10]{ACDS}).
\end{rem}

\section{Proof of Theorem~\ref{MainTheoremCor}}\label{sec:proof}

We first recall the notion of a divided difference.
For a function $f: \Rl \to \Cx$ the divided difference of the zeroth order~$f^{[0]}$ is
the function~$f$ itself. Let~$\lambda_0, \lambda_1, \ldots \in \Rl$
and let~$f \in C^n$. The divided difference $f^{[n]}$ of order~$n$ is
defined recursively by
\begin{align*}
  f^{[n]} \left( \lambda_0, \lambda_1, \tilde \lambda \right) =
  \begin{cases}\frac
    { f^{[n-1]} (\lambda_0, \tilde \lambda) - f^{[n-1]}(\lambda_1,
      \tilde \lambda)}{\lambda_0 - \lambda_1}, & \text{if}~\lambda_0
      \neq \lambda_1, \\ \frac {d}{d\lambda_1} f^{[n-1]} (\lambda_1,
    \tilde \lambda), & \text{if}~\lambda_0=\lambda_1,
  \end{cases}
\end{align*}
where~$\tilde\lambda = \left(\lambda_2, \ldots, \lambda_n \right) \in
\Rl^{n-1}$.

\begin{lem}\label{new_rep} Let $f$ be a Schwartz function on $\Rl$.

1. The divided difference $f^{[n]}$ can be written in the following form:
 \begin{equation}\label{dd_dfn}
\begin{aligned}
  f^{[n]} (\lambda_0, \ldots, \lambda_n) &=
  \int_{\Rl_+^{n+1}} \exp\left(i \sum_{j=0}^n t_j\lambda_j\right)\, \left( \mathcal F f\right)\left(\sum_{j=0}^n t_j\right) dt_0 dt_1 \cdots dt_n\\
  +(-1)^{n+1}&\int_{\Rl_+^{n+1}} \exp\left(-i \sum_{j=0}^n t_j\lambda_j\right)\, \left( \mathcal F f\right)\left(-\sum_{j=0}^n t_j\right) dt_0 dt_1 \cdots dt_n,
\end{aligned}
\end{equation}
where $\mathcal F f$ is the Fourier transform of the function $f$.

2. $f^{[n]}$ belongs to $\cC_n+ \cC_n$.

3. $f^{[n]}$ satisfies the assumptions of Theorem~\ref{qn_conv}.
\end{lem}

\begin{proof}
1. For every Schwartz function $f$ by~\cite[Lemma 2.3]{ACDS} we have the following representation
$$f^{[n]} (\lambda_0, \ldots, \lambda_n) = $$
$$=\int_{\Omega^n} e^{i (s_0-s_{1})\lambda_0}e^{i (s_1-s_{2})\lambda_1}\cdots
  e^{i (s_{n-1}-s_{n})\lambda_{n-1}}e^{i s_n\lambda_n}\, \left( \mathcal F f\right)(s_0) ds_0 ds_1 \cdots ds_n,$$
  where $\Omega^n = \{(s_0, s_1, \dots, s_n) \ : \ |s_n| \le \dots |s_1| \le |s_0|, {\rm sign} (s_n) = \cdots = {\rm sign} (s_0) \}.$

First, we write $\Omega^n$ as a union of
  $\Omega^n_+ = \{(s_0, s_1, \dots, s_n) \ : \ 0 \le s_n \le \dots s_1 \le s_0 \}$ and
  $\Omega^n_- := -\Omega^n_+$.
Second, we make a substitution $t_n=s_n$, $t_k=s_k-s_{k+1}$, $0\le k <n$ in the above integral. We obtain
\begin{align*}
  f^{[n]} (\lambda_0, \ldots, \lambda_n) &=
  \int_{\Rl_+^{n+1}} \exp\left(i \sum_{j=0}^n t_j\lambda_j\right)\, \left( \mathcal F f\right)\left(\sum_{j=0}^n t_j\right) dt_0 dt_1 \cdots dt_n\\
  &+\int_{\Rl_-^{n+1}} \exp\left(i \sum_{j=0}^n t_j\lambda_j\right)\, \left( \mathcal F f\right)\left(\sum_{j=0}^n t_j\right) dt_0 dt_1 \cdots dt_n,
\end{align*}
where $\Rl_-^{n+1}:= - \Rl_+^{n+1}$.

Next,
\begin{align*}
f^{[n]} (\lambda_0, \ldots, \lambda_n) =
  &\int_{\Rl_+^{n+1}} \exp\left(i \sum_{j=0}^n t_j\lambda_j\right)\, \left( \mathcal F f\right)\left(\sum_{j=0}^n t_j\right) dt_0 dt_1 \cdots dt_n\\
  +(-1)^{n+1}&\int_{\Rl_+^{n+1}} \exp\left(-i \sum_{j=0}^n t_j\lambda_j\right)\, \left( \mathcal F f\right)\left(-\sum_{j=0}^n t_j\right) dt_0 dt_1 \cdots dt_n.
\end{align*}
which proves the first assertion.

2. According to the latter formula for every Schwartz function $f$ the corresponding functions $a_j$ (from the representation~\eqref{Crep}) are $e^{\pm i t_j\lambda_j}$,
so they are bounded and continuous on $\Rl$, $j=0,...n$.
The corresponding measure $\nu$ on $\Rl_+^{n+1}$ is such that
 $d\nu(t) = \left( \mathcal F f\right)(t_0+\dots +t_n)dt_0 dt_1 \cdots dt_n$ for the first integral and
 $d\nu(t) = (-1)^{n+1}\left( \mathcal F f\right)(-t_0-\dots -t_n)dt_0 dt_1 \cdots dt_n$ for the second one.

 Since the Fourier transform of a Schwartz function is a Schwartz function itself, it follows that
 the measure $|\nu|$ is finite for both integrals and the condition~\eqref{phi_norm} is satisfied. Hence, $f^{[n]} \in \cC_n+ \cC_n$.

3. As was shown in the first and the second parts of the proof the measure $\nu$ (corresponding to $f^{[n]}$ in representation~\eqref{Crep})
is absolutely continuous with either $\nu'(t) = \left( \mathcal F f\right)(t_0+\dots +t_n)$ or $\nu'(t) = (-1)^{n+1}\left( \mathcal F f\right)(-t_0-\dots -t_n)$.
As was explained above, $\mathcal F f$ is a Schwartz function. Hence, $\nu'$ satisfies~\eqref{sw1}.

The functions $(t_0, \dots, t_n) \mapsto e^{\pm i t_j H}$ are bounded on $\Rl^{n+1}$, $j=0,...,n$. Also, all of them are Lipschitz.
Indeed, for every $a,b \in \Rl$ we have
$$\|e^{i aH} - e^{i bH}\|_\infty = \|1 - e^{i (b-a)H}\|_\infty  \le \sup_{|x| \le \|H\|_\infty |b-a|} |1 - e^{i x}| \le \|H\|_\infty |b-a|.$$

This completes the verification of the conditions of Theorem 3.8 and the proof of Lemma 4.1.

\end{proof}

The following result is a straightforward corollary of Lemma~\ref{new_rep} and Corollary~\ref{int_cor}.
\begin{cor}\label{div_diff}
Let $1\leq p_j \leq \infty$, $1\leq j\leq n$ and let $p$ be such that $\frac 1{p} = \sum_{j=1}^n \frac 1{p_j}$.
Let $x_j \in \mL_{p_j,\infty}$ and let~$H\in B(\sH)$.
For every Schwartz function $f$
the integral $T_{f^{[n]}}$ exists and the element
 $T_{f^{[n]}}(x_1, \dots, x_n)$ (which is given by by the integral in~\eqref{new_int}) is the limit of integral sums of the form~\eqref{S_m} with respect to the quasi-norm of $\mL_{p,\infty}$.
\end{cor}

%
%

Before we proceed with the proof of Theorem~\ref{M2} we state the following technical result, which is used below.

\begin{lem}\label{ph}
Let $f$ be a Schwartz function and set
$$\phi(\lambda_0, \lambda_1, \ldots, \lambda_{n-1})=-i  f^{[n]} \left( \lambda_0, \lambda_0, \lambda_1, \ldots,\lambda_{n-1} \right).$$
The following representations of the function $\phi$ hold:

1. \begin{align*}
\phi(\lambda_0, \lambda_1, \ldots, \lambda_{n-1})&=
  \int_{\Rl_+^{n}} t_0 \exp\left(i \sum_{j=0}^{n-1} t_j\lambda_j\right)\, \left( \mathcal F f\right)\left(\sum_{j=0}^{n-1} t_j\right) dt_0 dt_1 \cdots dt_{n-1}\\
  +(-1)^{n+1}&\int_{\Rl_+^{n}} t_0 \exp\left(-i \sum_{j=0}^{n-1} t_j\lambda_j\right)\, \left( \mathcal F f\right)\left(-\sum_{j=0}^{n-1} t_j\right) dt_0 dt_1 \cdots dt_{n-1}.
\end{align*}
In particular, $\phi$ satisfies the assumptions of Theorem~\ref{qn_conv}.

2. $$\phi(\lambda_0, \lambda_1, \ldots, \lambda_{n-1}) = \int_{\mathbb S^{n-1}} s_0 f^{(n)}\left(\sum_{j=0}^{n-1} \lambda_j s_j \right) d\sigma_{n-1},$$
where
$$\mathbb S^n = \left\{ \left( s_0, \ldots, s_n \right) \in \Rl^{n + 1}_+\ : \  \sum_{j = 0}^n s_j = 1\right\}$$
and $d\sigma_n$ is a finite measure on $\mathbb S^n$ defined by requiring that for every continuous function~$g: \Rl^{n+1} \mapsto \Cx$
the following equality holds:
\begin{equation*}
\label{minv}
\int_{S^n} g(s_0, \ldots, s_n)\, d\sigma_n =
\int_{R^n} g\left( s_0, \ldots, s_{n-1}, 1 - \sum_{j =
    0}^{n-1} s_j \right)\, ds_0 ds_1 \cdots ds_{n-1},
\end{equation*}
where
\begin{equation*}
  \label{Rn}
  R^n = \left\{(s_0, \ldots, s_{n-1})
    \in \Rl^n:\ \ \sum_{j = 0}^{n-1} s_j \leq 1,\ \ s_j \geq 0,\ 0
    \leq j \leq n\right\}.
\end{equation*}
\end{lem}

\begin{proof}
1. By the definition of the divided difference and Lemma~\ref{new_rep}, we have
 \begin{align*}
 \phi(\lambda_0, \lambda_1, \ldots, \lambda_{n-1}) &= -i f^{[n]} \left( \lambda_0, \lambda_0, \lambda_1, \ldots,\lambda_{n-1} \right)\\
 &=-i \frac{\partial}{\partial \lambda_0} f^{[n-1]} \left( \lambda_0, \lambda_1, \ldots,\lambda_{n-1} \right)\\
 =-i \frac{\partial}{\partial \lambda_0}\int_{\Rl_+^{n}} &\exp\left(i \sum_{j=0}^{n-1} t_j\lambda_j\right)\, \left( \mathcal F f\right)\left(\sum_{j=0}^{n-1} t_j\right) dt_0 dt_1 \cdots dt_{n-1}\\
  -i (-1)^{n}\frac{\partial}{\partial \lambda_0}\int_{\Rl_+^{n}} &\exp\left(-i \sum_{j=0}^{n-1} t_j\lambda_j\right)\, \left( \mathcal F f\right)\left(-\sum_{j=0}^{n-1} t_j\right) dt_0 dt_1 \cdots dt_{n-1}\\
=\int_{\Rl_+^{n}} &t_0 \exp\left(i \sum_{j=0}^{n-1} t_j\lambda_j\right)\, \left( \mathcal F f\right)\left(\sum_{j=0}^{n-1} t_j\right) dt_0 dt_1 \cdots dt_{n-1}\\
  +(-1)^{n+1}\int_{\Rl_+^{n}} &t_0 \exp\left(-i \sum_{j=0}^{n-1} t_j\lambda_j\right)\, \left( \mathcal F f\right)\left(-\sum_{j=0}^{n-1} t_j\right) dt_0 dt_1 \cdots dt_{n-1}.
\end{align*}
The arguments similar to that of Lemma~\ref{new_rep}(2) prove that $\phi$ satisfies the assumptions of Theorem~\ref{qn_conv}.

2. By~\cite[Lemma 5.1]{PSS} for every $f \in C^n$ we have the following representation:
$$f^{[n]}(\lambda_0, \lambda_1, \ldots, \lambda_{n}, ) = \int_{\mathbb S^{n}} f^{(n)}\left(\sum_{j=0}^{n} \lambda_j t_j \right) d\sigma_{n}.$$
Due to the latter formula (see also~\cite[Chapter IV, \S 7(a)]{DVL}) we obtain
$$f^{[n]}(\lambda_0, \lambda_1, \ldots, \lambda_{n} ) =f^{[n]}(\lambda_1, \ldots, \lambda_{n}, \lambda_0).$$
Hence, by the first part of this lemma and the definition of $d\sigma_{n}$, we have
\begin{align*}
 \phi(\lambda_0, \lambda_1, \ldots, \lambda_{n-1}) = -i &f^{[n]} \left( \lambda_0, \lambda_1, \ldots,\lambda_{n-1}, \lambda_0 \right)\\
 = -i \int_{\mathbb S^{n}} &f^{(n)}\left(\lambda_0(t_0+t_n) + \sum_{j=1}^{n-1} \lambda_j t_j \right) d\sigma_{n}\\
 = -i \int_{R^{n}} &f^{(n)}\left(\lambda_0(t_0+1- \sum_{j=0}^{n-1}  t_j) + \sum_{j=1}^{n-1} \lambda_j t_j \right) dt_0dt_1 \dots dt_{n-1}\\
 = -i \int_{R^{n}} &f^{(n)}\left(\lambda_0(1- \sum_{j=1}^{n-1}  t_j) + \sum_{j=1}^{n-1} \lambda_j t_j \right) dt_0dt_1 \dots dt_{n-1}.
\end{align*}

By the definition of $R^n$ we obtain that $0\le t_0 \le 1- \sum_{j=1}^{n-1}  t_j$ and $\sum_{j=1}^{n-1} t_j \le 1$.
So $(t_1, t_2, \dots, t_{n-1}) \in R^{n-1}$ and we have
\begin{align*}
 &\phi(\lambda_0, \lambda_1, \ldots, \lambda_{n-1})\\
 &= -i \int_{R^{n-1}} \int_0^{1- \sum_{j=1}^{n-1}  t_j} dt_0 f^{(n)}\left(\lambda_0(1- \sum_{j=1}^{n-1}  t_j) + \sum_{j=1}^{n-1} \lambda_j t_j \right) dt_1 \dots dt_{n-1}\\
 &= -i \int_{R^{n-1}} \left(1- \sum_{j=1}^{n-1}  t_j \right)f^{(n)}\left(\lambda_0(1- \sum_{j=1}^{n-1}  t_j) + \sum_{j=1}^{n-1} \lambda_j t_j \right) dt_1 \dots dt_{n-1}.
\end{align*}
Next, we make the following substitution: $s_0=1- \sum_{j=1}^{n-1}  t_j$, $s_k=t_k$, $1\le k\le n-2$.
Note that, $t_{n-1}=1- \sum_{j=0}^{n-2}  s_j$. Also note that, $\sum_{j=0}^{n-2}  s_j \le 1$, so $(s_0,\dots, s_{n-2})\in R^{n-1}$. The Jacobian of this substitution is $(-1)^{n+1}$ and so
$ds_0 \dots ds_{n-2}= dt_1 \dots dt_{n-1}$.
Hence,
\begin{align*}
 &\phi(\lambda_0, \lambda_1, \ldots, \lambda_{n-1})\\
 &= -i \int_{R^{n-1}} s_0 f^{(n)}\left(\lambda_0s_0 + \sum_{j=1}^{n-2} \lambda_j s_j  + \lambda_{n-1} (1- \sum_{j=0}^{n-2}  s_j)\right) ds_0 \dots ds_{n-2}\\
 &= -i \int_{\mathbb S^{n-1}} s_0 f^{(n)}\left(\sum_{j=0}^{n-1} \lambda_j s_j \right) d\sigma_{n-1},
 \end{align*}
 by the definition of the measure $d\sigma_{n}$.

\end{proof}

Finally, we are able to present the proof of Theorem~\ref{MainTheoremCor} stated in Preliminaries.

\begin{proof} [Proof of Theorem~\ref{MainTheoremCor}]
 Let $f : \Rl \to \Rl$ be a Schwartz function and $H \in B(\sH)$, $V \in \mL_{n,\infty}$ be self-adjoint operators.
For $p>n$ we have $\mL_{n,\infty} \subset \mL_p$ and, so $V \in \mL_p$. By~\cite[Theorem~5.7]{ACDS}
and Lemma~\ref{GF} the function $t \mapsto f(H+tV)$
is $n$-times differentiable at $H$ and
\begin{equation*}
    \frac {d^n}{dt^n} \left[ f(H+tV) \right]\biggr|_{t = 0} = n!\, \hat T_{f^{[n]}} \bigl(
    \underbrace{V, \ldots, V}_{\text{$n$-times}} \bigr),
\end{equation*}
where $\hat T_{f^{[n]}} : \mL_{p}^n \to \mL_{p/n}$ is a multiple operator integral in the sense of Definition~\ref{MOI}.
Here by $\mL_{p}^n$ we denote $\underbrace{\mL_{p} \times \cdots \times \mL_{p}}_{\text{$n$-times}}.$

We now show that if $H \in B(\sH)$ and $V \in \mL_{n,\infty}$ then
$$\hat T_{f^{[n]}}\bigl(\underbrace{V, \ldots, V}_{\text{$n$-times}} \bigr) =
T_{f^{[n]}} \bigl(\underbrace{V, \ldots, V}_{\text{$n$-times}} \bigr),$$
where $T_{f^{[n]}}: \mL_{n,\infty}^n \to \mL_{1,\infty}$ is a multiple operator integral in the sense of Definition~\ref{quasiMOI}.
Indeed, we have
$$\|\hat T_{f^{[n]}}\bigl(\underbrace{V, \ldots, V}_{\text{$n$-times}} \bigr) - T_{f^{[n]}}\bigl(\underbrace{V, \ldots, V}_{\text{$n$-times}} \bigr)\|_{\mL_{p/n}} \le$$
$$\le \|\hat T_{f^{[n]}}\bigl(\underbrace{V, \ldots, V}_{\text{$n$-times}} \bigr) - S_m\|_{\mL_{p/n}} +\|T_{f^{[n]}}\bigl(\underbrace{V, \ldots, V}_{\text{$n$-times}} \bigr)- S_m\|_{\mL_{1,\infty}},$$
where $S_m$ are the integral sums of the form~\eqref{S_m} corresponding to the function 
$$ (t_0, t_1, \dots, t_n) \mapsto e^{i t_0 H}V e^{i t_1 H}V \cdots e^{i t_{n-1} H} V  e^{i t_n H} \left( \mathcal F f \right)(t_0+\dots+t_n).$$
Now, the first term tends to zero by Theorem~\ref{bohner} and the second term tends to zero by Corollary~\ref{div_diff}.

Therefore,
\begin{equation}
    \label{SSFFrechetTemp}
    \frac {d^n}{dt^n} \left[ f(H+tV) \right]\biggr|_{t = 0} = n!\, T_{f^{[n]}} \bigl(
    \underbrace{V, \ldots, V}_{\text{$n$-times}} \bigr),
\end{equation}
where $T_{f^{[n]}} : \mL_{n,\infty}^n \to \mL_{1,\infty}$ is a multiple operator integral in the sense of Definition~\ref{quasiMOI}.
Hence, $\frac{d^n}{dt^{n}}\left[f\left(H+tV\right)\right]\biggr|_{t = 0}$ belongs to $\mL_{1,\infty}$.

Next, for every trace $\tau$ on $\mL_{1,\infty}$ we obtain
\begin{equation}\label{e1}
\frac 1{n!} \, \tau \left(\frac{d^n}{dt^n}\left[f\left(H+tV\right)\right]\biggr|_{t = 0} \right)
=\tau\left(T_{f^{[n]}}\bigl(\underbrace{V,\ldots,V}_{\text{$n$-times}}\bigr)\right).
\end{equation}

Combining Definition~\ref{quasiMOI} with Lemma~\ref{new_rep} yields
\begin{equation}\label{e0}
\begin{aligned}
 T_{f^{[n]}} \bigl(\underbrace{V, \ldots, V}_{\text{$n$-times}}\bigr)
 =&\int_{\Rl_+^{n+1}} e^{i t_0 H}V e^{i t_1 H}V \cdots e^{i t_{n-1} H} V e^{i t_n H} \left( \mathcal F f \right)(\sum_{k=0}^nt_k)\,\prod_{k=0}^ndt_k\\
 +(-1)^{n+1}\int_{\Rl_+^{n+1}}& e^{-i t_0 H}V e^{-i t_1 H}V \cdots e^{-i t_{n-1} H} V e^{-i t_n H} \left( \mathcal F f \right)(-\sum_{k=0}^nt_k)\,\prod_{k=0}^ndt_k.
\end{aligned}
\end{equation}

Recall that all integrals are understood in the sense of Definition~\ref{op_int}.

We now consider the first integral from the latter expression. The second integral is treated similarly.
By Theorem~\ref{qn_conv} this integral is the limit of integral sums of the form~\eqref{S_m}.
Hence, for every bounded trace $\tau$ on $\mL_{1,\infty}$ by Theorem~\ref{swap} we have
\begin{equation}\label{eq0}
\begin{aligned}
&\tau \left( \int_{\Rl_+^{n+1}} e^{i t_0 H}V e^{i t_1 H}V \cdots e^{i t_{n-1} H} V e^{i t_n H} \left( \mathcal F f \right)(t_0+\dots+t_n)\, dt_0 dt_1 \cdots dt_n \right) \\
&=\int_{\Rl_+^{n+1}} \tau \left( e^{i t_0 H}V e^{i t_1 H}V \cdots e^{i t_{n-1} H} V e^{i t_n H}  \right)\left( \mathcal F f \right)(t_0+\dots+t_n)\, dt_0 dt_1 \cdots dt_n\\
&=\int_{\Rl_+^{n+1}} \tau \left( e^{i (t_0+t_n) H}V e^{i t_1 H}V \cdots e^{i t_{n-1} H} V   \right)\left( \mathcal F f \right)(t_0+\dots+t_n)\, dt_0 dt_1 \cdots dt_n,
\end{aligned}
\end{equation}
where the latter equality is due to the following property of traces:
$\tau(AB) = \tau(BA)$ for all $A \in \mL_{1,\infty}$ and $B\in B(\sH)$.

Arguing as in the proof of Lemma~\ref{new_rep}(3) one can show that
the function
$$ (t_0, t_1, \dots, t_n) \mapsto e^{i (t_0+t_n) H}V e^{i t_1 H}V \cdots e^{i t_{n-1} H} V   \left( \mathcal F f \right)(t_0+\dots+t_n)$$
satisfies the assertions of Theorem~\ref{qn_conv}.
Therefore, it follows from Theorem~\ref{swap} that one can swap back the trace and the integral.

Hence,
\begin{equation}\label{eq00}
\begin{aligned}
&\tau \left( \int_{\Rl_+^{n+1}} e^{i t_0 H}V e^{i t_1 H}V \cdots e^{i t_{n-1} H} V e^{i t_n H} \left( \mathcal F f \right)(t_0+\dots+t_n)\, dt_0 dt_1 \cdots dt_n \right) \\
=&\tau \left(\int_{\Rl_+^{n+1}}  e^{i (t_0+t_n) H}V e^{i t_1 H}V \cdots e^{i t_{n-1} H} V   \left( \mathcal F f \right)(t_0+\dots+t_n)\, dt_0 dt_1 \cdots dt_n\right).
\end{aligned}
\end{equation}

Next, we claim that
\begin{equation}\label{eq000}
\begin{aligned}
 &\int_{\Rl_+^{n+1}}  e^{i (t_0+t_n) H}V e^{i t_1 H}V \cdots e^{i t_{n-1} H} V   \left( \mathcal F f \right)(t_0+\dots+t_n)\, dt_0 dt_1 \cdots dt_n\\
=&\left(\int_{\Rl_+^{n+1}}  e^{i (t_0+t_n) H}V e^{i t_1 H}V \cdots e^{i t_{n-1} H} \left( \mathcal F f \right)(t_0+\dots+t_n)\, dt_0 dt_1 \cdots dt_n \right)\cdot V,
\end{aligned}
\end{equation}
where the integral $I$ on the left-hand side is an integral in $\mL_{1,\infty}$ and the integral $J$ on the right-hand side is an integral in $\mL_{\frac{n}{n-1},\infty}$.
Indeed, if $S_m$ are the integral sums of $J$ of the form~\eqref{S_m}, then $S_m V$ are the integral sums of $I$ and
$$\|I-JV\|_{\mL_{1,\infty}} \le C\|I-S_mV\|_{\mL_{1,\infty}}+ C\|J-S_m\|_{\mL_{\frac{n}{n-1}}}\|V\|_{\mL_{n,\infty}}\to 0. $$
In the latter integral we make the following substitution $s_0=t_0+t_n$, $s_k=t_k$, $0<k\le n$.
Noting that $0\le t_n \le s_0$ (since $t_0 \ge0$) we obtain
\begin{equation}\label{eq01}
\begin{aligned}
&\int_{\Rl_+^{n+1}}  e^{i (t_0+t_n) H}V e^{i t_1 H}V \cdots e^{i t_{n-1} H} \left( \mathcal F f \right)(t_0+\dots+t_n)\, dt_0 dt_1 \cdots dt_n\\
=&\int_{\Rl_+^{n}}  e^{i s_0 H}V e^{i s_1 H}V \cdots e^{i s_{n-1} H} \left( \mathcal F f \right)(s_0+\dots+s_{n-1})\, ds_0 ds_1 \cdots ds_{n-1} \int_0^{s_0} ds_n\\
=&\int_{\Rl_+^{n}}  s_0 e^{i s_0 H}V e^{i s_1 H}V \cdots e^{i s_{n-1} H} \left( \mathcal F f \right)(s_0+\dots+s_{n-1})\, ds_0 ds_1 \cdots ds_{n-1}.
\end{aligned}
\end{equation}

Combining~\eqref{eq00},~\eqref{eq000},~\eqref{eq01} yields
\begin{equation}
\begin{aligned}
&\tau \left( \int_{\Rl_+^{n+1}} e^{i t_0 H}V \cdots e^{i t_{n-1} H} V e^{i t_n H} \left( \mathcal F f \right)(t_0+\dots+t_n)\, dt_0 dt_1 \cdots dt_n \right) \\
=&\tau \left(\int_{\Rl_+^{n}}  t_0 e^{i t_0 H}V \cdots e^{i t_{n-1} H} \left( \mathcal F f \right)(t_0+\dots+t_{n-1})\, dt_0 dt_1 \cdots dt_{n-1}  \cdot V\right).
\end{aligned}
\end{equation}

Similarly for the second integral from~\eqref{e0} we have
\begin{equation}
\begin{aligned}
&\tau \left( \int_{\Rl_+^{n+1}} e^{-i t_0 H}V \cdots e^{-i t_{n-1} H} V e^{-i t_n H} \left( \mathcal F f \right)(-t_0-\dots-t_n)\, dt_0 dt_1 \cdots dt_n \right) \\
=&\tau \left(\int_{\Rl_+^{n}}  t_0 e^{-i t_0 H}V \cdots e^{-i t_{n-1} H} \left( \mathcal F f \right)(-t_0-\dots-t_{n-1})\, dt_0 dt_1 \cdots dt_{n-1}  \cdot V\right).
\end{aligned}
\end{equation}

Finally we obtain the following representation
\begin{align*}
&\tau \left( T_{f^{[n]}} \bigl(\underbrace{V, \ldots, V}_{\text{$n$-times}}\bigr) \right)\\
&=\tau \left(\int_{\Rl_+^{n}}  t_0 e^{i t_0 H}V \cdots e^{i t_{n-1} H} \left( \mathcal F f \right)(t_0+\dots+t_{n-1})\,\prod_{k=0}^{n-1}dt_k\cdot V\right)\\
 &+(-1)^{n+1}\tau \left(\int_{\Rl_+^{n}}  t_0 e^{-i t_0 H}V \cdots e^{-i t_{n-1} H} \left( \mathcal F f \right)(-t_0-\dots-t_{n-1})\,\prod_{k=0}^{n-1}dt_k\cdot V\right).
  \end{align*}

Consider the function
$$\phi(\lambda_0, \lambda_1, \ldots, \lambda_{n-1})=-i  f^{[n]} \left( \lambda_0, \lambda_0, \lambda_1, \ldots,\lambda_{n-1} \right).$$
By Lemma~\ref{ph}(1) the function $\phi$ satisfies the assumptions of Theorem~\ref{qn_conv}. Hence, by Corollary~\ref{int_cor} the integral $T_{\phi}$ is well-defined.
Moreover, according to the representation of $\phi$ given in Lemma~\ref{ph}(1) and Definition~\ref{quasiMOI} we obtain
\begin{align*}
&T_{\phi}\bigl(\underbrace{V, \ldots, V}_{\text{$(n-1)$-times}}\bigr)
= \int_{\Rl_+^{n}}  t_0 e^{i t_0 H}V \cdots e^{i t_{n-1} H} \left( \mathcal F f \right)(t_0+\dots+t_{n-1})\, dt_0 dt_1 \cdots dt_{n-1}\\
&+(-1)^{n+1}\int_{\Rl_+^{n}}  t_0 e^{-i t_0 H}V \cdots e^{-i t_{n-1} H} \left( \mathcal F f \right)(-t_0-\dots-t_{n-1})\, dt_0 dt_1 \cdots dt_{n-1}.
\end{align*}
Consequently,
\begin{equation}\label{SSFRemainderTempI}
 \tau \left( T_{f^{[n]}} \bigl(\underbrace{V, \ldots, V}_{\text{$n$-times}}\bigr) \right) =
\tau \left(  T_{\phi}\bigl(\underbrace{V, \ldots, V}_{\text{$(n-1)$-times}}\bigr) V \right).
\end{equation}

  By the second part of Lemma~\ref{ph} we have
$$\phi(\lambda_0, \lambda_1, \ldots, \lambda_{n-1}) = \int_{\mathbb S^{n-1}} s_0 f^{(n)}\left(\sum_{j=0}^{n-1} \lambda_j s_j \right) d\sigma_{n-1}.$$
  Hence, the function $\phi$ satisfies the conditions of~\cite[Theorem 5.3]{PSS}.

Next, by~\cite[Theorem 5.3]{PSS} for every $1 < p_j < \infty$, $1 \leq j \leq n-1$ such that $ 0 < \frac1p =\frac 1 {p_1} + \ldots + \frac 1 {p_{n-1}} < 1$
the following estimate holds:
$$ \left\| T_{\phi}\right\|_{\mL_{p_1} \times \ldots \times \mL_{p_{n-1}} \to \mL_{p}} \leq   c_p\left\| f^{(n)}\right\|_{L_\infty}.$$

Note that the multilinear operator integral $T_{\phi} : \mL_{p_1} \times \ldots \times \mL_{p_{n-1}} \to \mL_{p}$ from~\cite{PSS} is defined in a way which differs from ours.
However, it is proved in~\cite[Lemma 3.5]{PSS}, that for $\phi \in \cC_{n-1}$ this definition coincides with Definition~\ref{MOI}, that is
$$T_\phi (x_1, x_2, \ldots, x_n)
    = \int_{\Rl_+^{n+1}} a_0(H, s)\, x_1 \, a_1 (H, s)\, x_2 \cdot \ldots\cdot x_n a_n (H, s)\, d \nu(s),$$
where the integral on the right-hand side is the Bochner integral of the $\mL_{p}$-valued function.
Moreover, by Lemma~\ref{bohner} the Bochner integral coincides with the integral in the sense of Definition~\ref{op_int}.
So, our multiple operator integral $T_{\phi}$ coincides with that of~\cite{PSS}.

Next, Theorem~\ref{Gen_int} below yields that the operator $T_{\phi}$
acts from $\mL_{p_1,\infty} \times \ldots \times \mL_{p_{n-1},\infty} \to \mL_{p,\infty}$ and
$$ \left\| T_{\phi} \right\|_{\mL_{p_1,\infty} \times \ldots \times \mL_{p_{n-1},\infty} \to \mL_{p,\infty}} \leq   c'_p \left\| f^{(n)}\right\|_{L_\infty},$$
for every $1 < p_j < \infty$, $1 \leq j \leq n-1$ such that $ 0 < \frac1p =\frac 1 {p_1} + \ldots + \frac 1 {p_{n-1}} < 1$.  

In the particular case when $p_j= n$ for all $1 \leq j \leq n-1$ we obtain
\begin{equation}
    \label{SSFRemainderKey}
    \left\| T_{\phi}\bigl(\underbrace{V, \ldots, V}_{\text{$(n-1)
          $-times}}\bigr)
    \right\|_{\mL_{\frac n{n-1}, \infty}} \leq c_n\, \left\| f^{(n)}\right\|_{L_\infty}\,
    \left\| V \right\|_{\mL_{n,\infty}}^{n-1},\ \ 0 \leq t \leq 1.
  \end{equation}
  Combining~(\ref{SSFRemainderKey}) with~(\ref{SSFRemainderTempI}) and~\eqref{e1}, yields
$$\left|\tau \left(\frac{d^n}{dt^n}\left[f\left(H+tV\right)\right]\biggr|_{t = 0} \right)\right|
\leq \, c_n\, \left\|
        f^{(n)} \right\|_{L_\infty} \, \left\| V\right\|_{\mL_{n,\infty}}^n.$$
\end{proof}

\section{The result for the Dixmier-Macaev ideal}\label{sec:DM}

Recall that $\aM^{(q)}_{1,\infty}$ denotes the $q$-convexification of $\aM_{1,\infty}$.

First note, that in the case of weak ideals the $q$-convexification of $\mL_{1,\infty}$ is $\mL_{q,\infty}$,
in symbols $\mL_{1,\infty}^{(q)} = \mL_{q,\infty}$. However, this is not true in the case of the Dixmier-Macaev ideal.
Indeed, it is proved in~\cite[Proposition 4.9]{CRSS} that $\aM_{q,\infty} \subsetneq \aM_{1,\infty}^{(q)}$ for every $q>1$.

In a way similar to that of~\cite[Section 4.3]{CRSS} it can be shown that the norm in $\aM_{1,\infty}^{(q)}$ can be written in the following form:
\begin{equation}\label{norm}
 \|A\|_{\aM_{1,\infty}^{(q)}} = \sup_{s>1} (s-1)^{\frac1q} \|A\|_{\mL_{sq}}.
\end{equation}

Now we introduce the notion of multiple operator integral on $\aM_{p,\infty}$.
The following definition is similar to Definition~\ref{quasiMOI}

\begin{dfn}\label{quasiMOI1}
 Let $1\leq p_j \leq \infty$, $1\leq j\leq n$ and let $p$ be such that $\frac 1{p} = \sum_{j=1}^n\frac 1{p_j}$.
 Let $\phi\in \cC_n$ admits the representation~\eqref{Crep} with a measure $\nu$ being absolute continuous.
 A function $u : \Rl_+^{n+1} \to \aM_{p,\infty}$ is defined as follows:
\begin{equation}
 u_\phi(s_0, s_1, \dots, s_n) := a_0(H,s)x_1 a_1(H,s)x_2 \cdots x_n a_n(H,s) \nu'(s),
\end{equation}
where the function $\nu' : \Rl_+^{n+1} \to \Cx$ is the Radon-Nikodym derivative of $\nu$.

 If the integral $\int_{\Rl_+^{n+1}} u_\phi(s) \ ds$
 exists in the sense of Definition~\ref{op_int},
 then  we define the operator
 $$T_{\phi} : \aM_{p_1,\infty} \times \dots \times \aM_{p_n,\infty} \to \aM_{p,\infty}$$
as follows:
\begin{equation}
 T_{\phi}(x_0, x_1, \dots, x_n) := \int_{\Rl_+^{n+1}} u_\phi(s) \ ds.
\end{equation}
\end{dfn}

In a way similar to that of Section~\ref{sec:MOI} it can be showed that for a wide class of functions $\phi \in \cC_n$
(more specifically those described in Theorem~\ref{qn_conv})
the operator integral $T_\phi$ exists and is the limit of integral sums with respect to the norm of $\aM_{1, \infty}^{(q)}$.
In particular, if $f$ is Schwartz, then the integrals $T_{f^{[n]}}$ and $T_\phi$
(for $\phi(\lambda_0, \lambda_1, \ldots, \lambda_{n-1})= -i f^{[n]} \left( \lambda_0, \lambda_0, \lambda_1, \ldots,\lambda_{n-1} \right)$)
exist.

The following result proves the key estimate~\eqref{SSFRemainderKey} in the case of the Dixmier-Macaev ideal.

\begin{thm}\label{DM_thm}
 Let $n \in \N$ and let $f$ be a Schwartz function.
For every $H \in B(\sH)$ and $V \in \aM_{1,\infty}^{(n)}$ we have
$$ \left\| T_{\phi}(\underbrace{V, \ldots, V}_{\text{$(n-1)$-times}}) \right\|_{\aM_{1,\infty}^{( \frac{n}{n-1} )}}
\leq \,   c_n\, \left\| f^{(n)} \right\|_{L_\infty} \|V\|_{\aM_{1,\infty}^{(n)}}^{n-1}, $$
  where the constant~$c_n$ depends only on $n$ and,
  $T_{\phi}: \underbrace{\aM_{1,\infty}^{(n)} \times \cdots \times \aM_{1,\infty}^{(n)}}_{\text{$(n-1)$-times}} \to \aM_{1,\infty}^{( \frac{n}{n-1} )}$ is a multiple operator integral associated
  with~$H$, $V$ and the function
  $$\phi(\lambda_0, \lambda_1, \ldots, \lambda_{n-1})=
  -i f^{[n]} \left( \lambda_0, \lambda_0, \lambda_1, \ldots,\lambda_{n-1}
  \right).$$
\end{thm}

\begin{proof}
 According to~\eqref{norm} we have
 $$\left\| T_{\phi}(\underbrace{V, \ldots, V}_{\text{$(n-1)$-times}}) \right\|_{\aM_{1,\infty}^{( \frac{n}{n-1} )}}
 = \sup_{s>1} (s-1)^{\frac{n-1}{n}} \|T_{\phi}(\underbrace{V, \ldots, V}_{\text{$(n-1)$-times}})\|_{\mL_{\frac{sn}{n-1}}}.$$

By~\cite[Theorem 2.1]{PSS}, we have
\begin{align*}
 \left\| T_{\phi}(\underbrace{V, \ldots, V}_{\text{$(n-1)$-times}}) \right\|_{\aM_{1,\infty}^{( \frac{n}{n-1} )}}
 &\le c\, \left\| f^{(n)} \right\|_{L_\infty} \sup_{s>1} (s-1)^{\frac{n-1}{n}} \|V\|_{\mL_{\frac{sn}{n-1}(n-1)}}^{n-1}\\
 &= c\, \left\| f^{(n)} \right\|_{L_\infty} \left(\sup_{s>1} (s-1)^{\frac1{n}} \|V\|_{\mL_{sn}}\right)^{n-1}\\
 &= c\, \left\| f^{(n)} \right\|_{L_\infty} \|V\|_{\aM_{1,\infty}^{(n)}}^{n-1}.
\end{align*}
\end{proof}

Using the latter result we are able to extend Theorems~\ref{MainTheoremCor} and~\ref{M2} to the Dixmier-Macaev ideal.

\begin{thm}\label{MR_DM}
  Let $n\in\N$. Let~$H \in B(\sH)$ be a self-adjoint operator
  and let~$V$ be a self-adjoint operator in $\aM_{1,\infty}^{(n)}$.
  For every Schwartz function $f$ the operator
  $\frac{d^n}{dt^n}\left[f\left(H+tV\right)\right]\biggr|_{t = 0}$
  belongs to $\aM_{1,\infty}$.
Moreover, there is a constant $c_n$ depending only on $n$ such that for every bounded trace $\tau$ on $\aM_{1,\infty}$ the estimate
\begin{equation}
\left|\tau \left(\frac{d^n}{dt^n}\left[f\left(H+tV\right)\right]\biggr|_{t = 0} \right) \right| \leq \, c_n\, \left\|
        f^{(n)} \right\|_{L_\infty} \, \left\| V\right\|_{\aM_{1,\infty}^{(n)}}^n
\end{equation} holds.
\end{thm}

The proof is a verbatim repetition of that of Theorem~\ref{MainTheoremCor}, with the only difference that we use Theorem~\ref{DM_thm}
instead of the interpolation argument used in the proof of Theorem~\ref{MainTheoremCor}.

Finally, using Theorem~\ref{MR_DM} we can prove Theorem~\ref{M1}. 
The proof is similar to that of Theorem~\ref{M2} and therefore omitted.

%

\section{The case of an unbounded operator $H$}\label{sec:Unbounded}

In the present section we prove the perturbation formulae for the unbounded operator $H$.
We deal with the class of functions $f$ considered by Koplienko~\cite{Koplienko},
that is the class of rational functions with non-real poles which
are bounded at infinity. Note that every function from this class belongs to the span of the following set:
$$\left\{ \lambda \mapsto (z-\lambda)^{-m} \ : \ m \in \N, z\notin \Rl\right\}.$$

We start with the representation of the divided difference for this class of functions.
\begin{lem}\label{6.1}
 Let $m \in \N$ and $z\notin \Rl$.
For the function $f : \Rl \to \Cx$ given by $f(\lambda) = (z-\lambda)^{-m}$ the $n^{\mathrm{th}}$ divided difference of $f$ can be written in the following form:
$$f^{[n]}(\lambda_0, \lambda_1, \ldots, \lambda_{n-1}, \lambda_n)
=\sum_{\genfrac{}{}{0pt}{}{1 \le m_0,\dots, m_n \le m}{m_0+\cdots +m_n = m+n}} \prod_{i=0}^{n}(z-\lambda_i)^{-m_i}.$$
\end{lem}

\begin{proof}
 We prove the formula by induction. For $n=0$ the formula is evidently correct.
 Assume that
 $$f^{[n-1]}(\lambda_0, \lambda_1, \ldots, \lambda_{n-1})
=\sum_{\genfrac{}{}{0pt}{}{1 \le m_0,\dots, m_{n-1} \le m}{m_0+\cdots +m_{n-1} = m+n-1}} \prod_{i=0}^{n-1}(z-\lambda_i)^{-m_i}.$$

We have
\begin{align*}
 &f^{[n]}(\lambda_0, \lambda_1, \ldots, \lambda_{n-1}, \lambda_n) \\
 &= \frac{f^{[n-1]}(\lambda_0, \lambda_2, \ldots, \lambda_{n-1}, \lambda_n) - f^{[n-1]}(\lambda_1, \lambda_2, \ldots, \lambda_{n-1}, \lambda_n)}{\lambda_0 - \lambda_1}\\
 &=\frac1{\lambda_0 - \lambda_1} \left( \sum_{\genfrac{}{}{0pt}{}{1 \le m_0,\dots, m_{n-1} \le m}{m_0+\cdots +m_{n-1} = m+n-1}} (z-\lambda_0)^{-m_0}\prod_{i=1}^{n-1}(z-\lambda_{i+1})^{-m_i}\right.\\
 &\qquad - \left.\sum_{\genfrac{}{}{0pt}{}{1 \le m_0,\dots, m_{n-1} \le m}{m_0+\cdots +m_{n-1} = m+n-1}} (z-\lambda_1)^{-m_0}\prod_{i=1}^{n-1}(z-\lambda_{i+1})^{-m_i}\right)\\
 &=\frac1{\lambda_0 - \lambda_1} \sum_{\genfrac{}{}{0pt}{}{1 \le m_0,\dots, m_{n-1} \le m}{m_0+\cdots +m_{n-1} = m+n-1}} 
 \left( (z-\lambda_0)^{-m_0} - (z-\lambda_1)^{-m_0}\right)\prod_{i=1}^{n-1}(z-\lambda_{i+1})^{-m_i}.
\end{align*}
Note that,
$$\frac{(z-\lambda_0)^{-m_0} - (z-\lambda_1)^{-m_0}}{\lambda_0 - \lambda_1} = \sum_{\genfrac{}{}{0pt}{}{1\le k,l \le m_0}{k+l =m_0+1}}(z-\lambda_0)^{-k} (z-\lambda_1)^{-l}.$$
Therefore,
\begin{align*}
 &f^{[n]}(\lambda_0, \lambda_1, \ldots, \lambda_{n-1}, \lambda_n) \\
 &= \sum_{\genfrac{}{}{0pt}{}{1 \le m_0,\dots, m_{n-1} \le m}{m_0+\cdots +m_{n-1} = m+n-1}} 
 \sum_{\genfrac{}{}{0pt}{}{1\le k,l \le m_0}{k+l =m_0+1}}(z-\lambda_0)^{-k} (z-\lambda_1)^{-l}\prod_{i=1}^{n-1}(z-\lambda_{i+1})^{-m_i}.
\end{align*}
Rename the variables as follows: $m_0:=k, m_1:=l, m_i:=m_{i-1}$, $1\le i \le n$. Finally, we obtain
$$ f^{[n]}(\lambda_0, \lambda_1, \ldots, \lambda_{n-1}, \lambda_n)
=\sum_{\genfrac{}{}{0pt}{}{1 \le m_0,\dots, m_{n} \le m}{m_0+\cdots +m_n = m+n}} \prod_{i=0}^{n}(z-\lambda_i)^{-m_i}.$$
\end{proof}


If $f$ is as in the previous lemma, then for every $0\leq p_i \leq \infty$ and $x_i \in \mL_{p_i,\infty}$, $1\leq i\leq n$ we have
$$T_{f^{[n]}}(x_1, \dots, x_n) = \sum_{\genfrac{}{}{0pt}{}{1 \le m_0,\dots, m_{n} \le m}{m_0+\cdots +m_n = m+n}} (z-\lambda_0)^{-m_0}\prod_{i=1}^{n}x_i(z-\lambda_i)^{-m_i},$$
where $T_{f^{[n]}}: \mL_{p_1,\infty} \times \dots \times \mL_{p_n,\infty} \to \mL_{p,\infty}$ is the operator integral in the sense of Definition~\ref{quasiMOI}.

The following theorem is a cornerstone estimate in the proof of Theorem~\ref{M4}.

\begin{thm}\label{unb}
 Let $f$ be a rational function with non-real poles which
is bounded at infinity. If $H, V$ are self-adjoint operators
such that~$V \in \mL_{n,\infty}$, then

1. The operator $\frac {d^n}{dt^n} \left[ f\left( H+tV\right) \right] \biggr|_{t = 0}$ belongs to $\mL_{1,\infty}$;

2. For every bounded trace $\tau$ on $\mL_{1,\infty}$ we have
$$\left| \tau\left( \frac {d^n}{dt^n} \left[ f\left( H+tV\right) \right] \biggr|_{t = 0} \right) \right| \le c_n\, \left\|
        f^{(n)} \right\|_{L_\infty} \, \left\| V\right\|_{\mL_{n,\infty}}^n.$$
\end{thm}

\begin{proof}
 It is sufficient to prove the first assertion for the function $g(\lambda) = (z-\lambda)^{-m}$.

 Due to~\cite[formula 2.4]{Koplienko} we have
$$\frac 1{n!} \, \frac {d^n}{dt^n} \left[ g\left( H+tV\right) \right] \biggr|_{t = 0} 
= \sum_{\genfrac{}{}{0pt}{}{1 \le m_0,\dots, m_{n} \le m}{m_0+\cdots +m_n = m+n}} (zI-H)^{-m_0} \prod_{i=1}^n V (zI-H)^{-m_i}.$$
Note that, every term in the latter sum is a product of bounded operators $(zI-H)^{-m_i}$ and $n$ operators $V$.
Since $V \in \mL_{n,\infty}$, it follows that the every term in this sum belongs to $\mL_{1,\infty}$,
that is $\frac {d^n}{dt^n} \left[ g\left( H+tV\right) \right] \biggr|_{t = 0} \in \mL_{1,\infty}$.
The first assertion has been proved.

Next, for every bounded trace $\tau$ on $\mL_{1,\infty}$ we have
\begin{align*}
 &\tau \left( \frac 1{n!} \, \frac {d^n}{dt^n} \left[ g\left( H+tV\right) \right] \biggr|_{t = 0} \right)\\
&= \tau \left( \sum_{\genfrac{}{}{0pt}{}{1 \le m_0,\dots, m_{n} \le m}{m_0+\cdots +m_n = m+n}} (zI-H)^{-m_0} \prod_{i=1}^n V (zI-H)^{-m_i}\right)\\
&=\tau \left( \sum_{\genfrac{}{}{0pt}{}{1 \le m_0,\dots, m_{n} \le m}{m_0+\cdots +m_n = m+n}} (zI-H)^{-m_0-m_n} \prod_{i=1}^{n-1} V (zI-H)^{-m_i} \cdot V\right),
\end{align*}
where the latter equality is due to the following property of traces:
$\tau(AB) = \tau(BA)$ for all $A \in \mL_{1,\infty}$ and $B\in B(\sH)$.

Note that
$$\sum_{\genfrac{}{}{0pt}{}{1 \le m_0,\dots, m_{n} \le m}{m_0+\cdots +m_n = m+n}} (zI-H)^{-m_0-m_n} \prod_{i=1}^{n-1} V (zI-H)^{-m_i}$$
equals to $T_\psi \bigl(\underbrace{V, \ldots, V}_{\text{$n-1$-times}}\bigr)$, where 
\begin{align*}
 \psi(\lambda_0, \lambda_1, \ldots, \lambda_{n-1})
&=\sum_{\genfrac{}{}{0pt}{}{1 \le m_0,\dots, m_{n} \le m}{m_0+\cdots +m_n = m+n}} (z-\lambda_0)^{-m_0-m_n} \prod_{i=1}^{n-1} (z-\lambda_i)^{-m_i}\\
&=g^{[n]}(\lambda_0, \lambda_1, \ldots, \lambda_{n-1}, \lambda_0),
\end{align*}

due to Lemma~\ref{6.1}.

Therefore,
\begin{equation}\label{eq006}
  \tau \left( \frac 1{n!} \, \frac {d^n}{dt^n} \left[ g\left( H+tV\right) \right] \biggr|_{t = 0} \right) =
  \tau \left( T_\psi \bigl(\underbrace{V, \ldots, V}_{\text{$n-1$-times}}\bigr) \cdot V\right)
\end{equation}
for every function $g \in\left\{ \lambda \mapsto (z-\lambda)^{-m} \ : \ m \in \N, z \notin \Rl\right\}.$
As was mentioned at the beginning of this section the function $f$ is a linear combination of such functions $g$. 
Therefore, due to the linearity of both sides of~\eqref{eq006}, we obtain
\begin{equation}\label{eq007}
  \tau \left( \frac 1{n!} \, \frac {d^n}{dt^n} \left[ f\left( H+tV\right) \right] \biggr|_{t = 0} \right) =
  \tau \left( T_\phi \bigl(\underbrace{V, \ldots, V}_{\text{$n-1$-times}}\bigr)  \cdot V\right),
\end{equation}
where $\psi(\lambda_0, \lambda_1, \ldots, \lambda_{n-1})=f^{[n]}(\lambda_0, \lambda_1, \ldots, \lambda_{n-1}, \lambda_0).$

Due to~\cite[Chapter IV, \S 7(a)]{DVL} we have
$$f^{[n]}(\lambda_0, \lambda_1, \ldots, \lambda_{n} ) =f^{[n]}(\lambda_1, \ldots, \lambda_{n}, \lambda_0).$$
So, $$\phi(\lambda_0, \lambda_1, \ldots, \lambda_{n-1})=f^{[n]}(\lambda_0, \lambda_0, \lambda_1, \ldots, \lambda_{n-1}).$$

By the second part of Lemma~\ref{ph} we have
$$\phi(\lambda_0, \lambda_1, \ldots, \lambda_{n-1}) = -i\int_{\mathbb S^{n-1}} s_0 f^{(n)}\left(\sum_{j=0}^{n-1} \lambda_j s_j \right) d\sigma_{n-1}.$$
  Hence, the function $\phi$ satisfies the conditions of~\cite[Theorem 5.3]{PSS}.

Next, by~\cite[Theorem 5.3]{PSS} for every $1 < p_j < \infty$, $1 \leq j \leq n-1$ such that $ 0 < \frac1p =\frac 1 {p_1} + \ldots + \frac 1 {p_{n-1}} < 1$
the following estimate holds:
$$ \left\| T_{\phi}\right\|_{\mL_{p_1} \times \ldots \times \mL_{p_{n-1}} \to \mL_{p}} \leq   c_p\left\| f^{(n)}\right\|_{L_\infty}.$$


Next, Theorem~\ref{Gen_int} below yields that the operator $T_{\phi}$
acts from $\mL_{p_1,\infty} \times \ldots \times \mL_{p_{n-1},\infty} \to \mL_{p,\infty}$ and
$$ \left\| T_{\phi} \right\|_{\mL_{p_1,\infty} \times \ldots \times \mL_{p_{n-1},\infty} \to \mL_{p,\infty}} \leq   c'_p \left\| f^{(n)}\right\|_{L_\infty},$$
for every $1 < p_j < \infty$, $1 \leq j \leq n-1$ such that $ 0 < \frac1p =\frac 1 {p_1} + \ldots + \frac 1 {p_{n-1}} < 1$.  

In the particular case when $p_j= n$ for all $1 \leq j \leq n-1$ we obtain
\begin{equation}
    \label{SSFRemainderKey0}
    \left\| T_{\phi}\bigl(\underbrace{V, \ldots, V}_{\text{$(n-1)
          $-times}}\bigr)
    \right\|_{\mL_{\frac n{n-1}, \infty}} \leq c_n\, \left\| f^{(n)}\right\|_{L_\infty}\,
    \left\| V \right\|_{\mL_{n,\infty}}^{n-1},\ \ 0 \leq t \leq 1.
  \end{equation}
  Combining~(\ref{SSFRemainderKey0}) with~~\eqref{eq007}, yields
$$\left|\tau \left(\frac{d^n}{dt^n}\left[f\left(H+tV\right)\right]\biggr|_{t = 0} \right)\right|
\leq \, c_n\, \left\|
        f^{(n)} \right\|_{L_\infty} \, \left\| V\right\|_{\mL_{n,\infty}}^n.$$
\end{proof}

Now using Theorem~\ref{unb} we can prove Theorem~\ref{M4}. 
The proof is similar to that of Theorem~\ref{M2} and therefore omitted.

In order to prove Theorem~\ref{M3} one needs a verbatim repetition of all the construction described in this section for the Dixmier-Macaev ideal $\aM_{1,\infty}$.
We omit futher details.

\section{Appendix (Polylinear interpolation)}\label{sec:App}

In this section we prove the result concerning the polylinear interpolation that we used in the preceding sections.

Let $0< p < \infty$, $0< q <\infty$. Define the Lorentz ideal $\mL_{pq}$ by setting (see e.g~\cite{LT2} and~\cite{DDP})
$$\mL_{p,q} := \left\{A \in B(\sH) \text{ is compact}:  \sum_{k=0}^\infty k^{\frac{q}{p}-1}\mu^q(k,A) < \infty\right\}.$$
Note that $\mL_{p,p} = \mL_p$.

We use the following
noncommutative version of the general Marcinkiewicz interpolation theorem for $\mL_{p, q}$ below.
It is a combination of~\cite[Theorem 5.3.2]{BL} and~\cite[Theorem 4.8]{Dirksen}.

\begin{thm}\label{interpolation}
Let $0< \alpha'\neq \alpha'' < \infty$, $0<\beta' \neq \beta'' < \infty$, $0 < \gamma', \gamma'', \delta', \delta'' \le \infty$ and $0 <\theta <1$.
Suppose that a linear operator $W$ acts from $\mL_{\alpha',\gamma'}$ to $\mL_{\beta',\delta'}$ and from $\mL_{\alpha'',\gamma''}$ to $\mL_{\beta'',\delta''}.$

If $\frac1\alpha =\frac{\theta}{\alpha'} +\frac{1-\theta}{\alpha''}$ and $\frac1\beta =\frac{\theta}{\beta'} +\frac{1-\theta}{\beta''}$,
then $W$ acts from $\mL_{\alpha,\gamma}$ to $\mL_{\beta,\gamma}$ for every $0<\gamma \le \infty$
and
$$\|W\|_{\mL_{\alpha,\gamma} \to \mL_{\beta,\gamma}} \le \|W\|_{\mL_{\alpha',\gamma'} \to \mL_{\beta',\delta'}}^\theta \|W\|_{\mL_{\alpha'',\gamma''} \to \mL_{\beta'',\delta''}}^{1-\theta}.$$
\end{thm}

The following theorem is the main result of the present section.

\begin{thm}\label{Gen_int}
If a multilinear operator~$R$ acts from $\mL_{\alpha_1} \times \ldots \times \mL_{\alpha_n}$ to $\mL_{\alpha}$
for every $0 < \alpha_j < \infty$, $1 \leq j \leq n$ such that $\frac1\alpha =\frac 1 {\alpha_1} + \ldots + \frac 1 {\alpha_n}$ and
$$ \left\| R\right\|_{\mL_{\alpha_1} \times \ldots \times \mL_{\alpha_n} \to \mL_{\alpha}} \leq \,   c_\alpha, $$
then
the operator~$R$ acts from
 $\mL_{\alpha_1,\infty} \times \ldots \times \mL_{\alpha_n,\infty}$ to $\mL_{\alpha,\infty}$
 for every $0 < \alpha_j < \infty$, $1 \leq j \leq n$ such that $\frac1\alpha =\frac 1 {\alpha_1} + \ldots + \frac 1 {\alpha_n}$ and
$$ \left\| R \right\|_{\mL_{\alpha_1,\infty} \times \ldots \times \mL_{\alpha_n,\infty} \to \mL_{\alpha,\infty}} \leq \,   c'_\alpha.$$
\end{thm}

\begin{proof}
1. Fix operators $x_j \in \mL_{\alpha_j}$, $j = 2, 3, \ldots, n$ and consider the operator
$W_1 : \mL_{\alpha_1} \to \mL_{\alpha}$ given by the formula
$$W_1(x_1) := R(x_1, x_2, ..., x_n).$$

Let us fix $0<\theta <1$ and find $\alpha', \alpha'', \beta', \beta''$ such that $\alpha' < \alpha_1 < \alpha''$ and
\begin{align*}
\frac1{\alpha_1} &=\frac{\theta}{\alpha'} +\frac{1-\theta}{\alpha''},\\
\frac1{\beta'} &=\frac1{\alpha'} +\sum_{j=2}^n \frac1{\alpha_j},\\
\frac1{\beta''} &=\frac1{\alpha''} +\sum_{j=2}^n \frac1{\alpha_j}.
\end{align*}

Note that, $\frac{\theta}{\beta'} +\frac{1-\theta}{\beta''}=\frac1{\alpha}$.

By the assumption on operator $R$ we have
$$ W_1 : \mL_{\alpha'} \to \mL_{\beta'} \ \text{and} \ W_1 : \mL_{\alpha''} \to \mL_{\beta''}.$$

Hence, by Theorem~\ref{interpolation} (with $\gamma'=\alpha'$, $\delta'=\beta'$, $\gamma''=\alpha''$, $\delta''=\beta''$ and $\gamma=\infty$)
$$ W_1 : \mL_{\alpha_1,\infty} \to \mL_{\alpha,\infty}$$
and
$$\|W_1\|_{\mL_{\alpha_1,\infty} \to \mL_{\alpha,\infty}} \le \|W_1\|_{\mL_{\alpha'} \to \mL_{\beta'}}^\theta \|W_1\|_{\mL_{\alpha''} \to \mL_{\beta''}}^{1-\theta}.$$

Therefore, by the definition of operator $W_1$ we conclude that
the operator~$R$ acts from  $\mL_{\alpha_1,\infty} \times \mL_{\alpha_2} \times \ldots \times \mL_{\alpha_n}$ to $\mL_{\alpha,\infty}$
for every $0 < \alpha_j < \infty$, $1 \leq j \leq n$ such that $\frac1\alpha =\frac 1 {\alpha_1} + \ldots + \frac 1 {\alpha_n}$ and
$$\left\| R (x_1, x_2, ..., x_n)\right\|_{\mL_{\alpha,\infty}} \leq \,   C \|x_1\|_{\mL_{\alpha_1,\infty}} \prod_{j=2}^n \|x_j\|_{\mL_{\alpha_j}},$$
for every $x_1 \in \mL_{\alpha_1,\infty}$ and $x_j \in \mL_{\alpha_j}$, $j = 2, 3, \ldots, n$.

2. Fix the operators $x_1 \in \mL_{\alpha_1,\infty}$ $x_j \in \mL_{\alpha_j}$, $j = 3, 4, \ldots, n$ and consider the operator
$W_2 : \mL_{\alpha_2} \to \mL_{\alpha,\infty}$ given by the formula
$$W_2(x_2) := R(x_1, x_2, ..., x_n).$$
Similarly, we fix $0<\theta <1$ and find $\alpha', \alpha'', \beta', \beta''$ such that $\alpha' < \alpha_2 < \alpha''$ and
\begin{align*}
\frac1{\alpha_2} &=\frac{\theta}{\alpha'} +\frac{1-\theta}{\alpha''},\\
\frac1{\beta'} &=\frac1{\alpha'} +\frac 1 {\alpha_1} +\sum_{j=3}^n \frac1{\alpha_j},\\
\frac1{\beta''} &=\frac1{\alpha''} +\frac 1 {\alpha_1} +\sum_{j=3}^n \frac1{\alpha_j}.
\end{align*}

Note that, $\frac{\theta}{\beta'} +\frac{1-\theta}{\beta''}=\frac1{\alpha}$.

It was proved in the first part, that $R$ acts from  $\mL_{\alpha_1,\infty} \times \mL_{\alpha_2} \times \ldots \times \mL_{\alpha_n}$ to $\mL_{\alpha,\infty}$
for every $0 < \alpha_j < \infty$, $1 \leq j \leq n$ such that $\frac1\alpha =\frac 1 {\alpha_1} + \ldots + \frac 1 {\alpha_n}$.
Hence,
$$ W_2 : \mL_{\alpha'} \to \mL_{\beta',\infty} \ \text{and} \ W_2 : \mL_{\alpha''} \to \mL_{\beta'',\infty}.$$

Therefore, by Theorem~\ref{interpolation} (with $\gamma'=\alpha'$, $\gamma''=\alpha''$ and $\delta'=\delta''=\gamma=\infty$)
$$ W_2 : \mL_{\alpha_2,\infty} \to \mL_{\alpha,\infty}$$
and
$$\|W_2\|_{\mL_{\alpha_2,\infty} \to \mL_{\alpha,\infty}} \le \|W_2\|_{\mL_{\alpha'} \to \mL_{\beta', \infty}}^\theta \|W_2\|_{\mL_{\alpha''} \to \mL_{\beta'', \infty}}^{1-\theta}.$$
In other words,
$R$ acts from  $\mL_{\alpha_1,\infty} \times \mL_{\alpha_2, \infty} \times \mL_{\alpha_3} \times \ldots \times \mL_{\alpha_n}$ to $\mL_{\alpha,\infty}$
for every $0 < \alpha_j < \infty$, $1 \leq j \leq n$ such that $\frac1\alpha =\frac 1 {\alpha_1} + \ldots + \frac 1 {\alpha_n}$ and
$$\left\| R (x_1, x_2, ..., x_n)\right\|_{\mL_{\alpha,\infty}}
\leq \,   C\|x_1\|_{\mL_{\alpha_1,\infty}} \|x_2\|_{\mL_{\alpha_2,\infty}}\prod_{j=3}^n \|x_j\|_{\mL_{\alpha_j}},$$
for every $x_1 \in \mL_{\alpha_1,\infty}$, $x_2 \in \mL_{\alpha_2,\infty}$ and $x_j \in \mL_{\alpha_j}$, $j = 3, \ldots, n$.

Repeating this procedure $(n-2)$ more times we end up having
$$\left\| R (x_1, x_2, ..., x_n)\right\|_{\mL_{\alpha,\infty}} \leq \,   c'_\alpha \prod_{j=1}^n \|x_j\|_{\mL_{\alpha_j,\infty}},$$
for every $x_j \in \mL_{\alpha_j,\infty}$, $j = 1, 2, \ldots, n$, that is
$$ \left\| R \right\|_{\mL_{\alpha_1,\infty} \times \ldots \times \mL_{\alpha_n,\infty} \to \mL_{\alpha,\infty}} \leq \,   c'_\alpha. $$
\end{proof}


\end{document}